\documentclass[a4paper,10pt]{article}

\usepackage{amsthm}

\usepackage{amsmath}

\usepackage{amsfonts}

\usepackage{amssymb}

\usepackage{amsmath,amsthm,amscd,amssymb}

\usepackage{latexsym}

\usepackage{graphicx}

\usepackage[all]{xy}

\theoremstyle{plain}

\newtheorem{thm}{\textbf{Theorem}}

\newtheorem{lem}{\textbf{Lemma}}

\newtheorem{df}{\textbf{Definition}}

\newtheorem{cor}{\textbf{Corollary}}

\newtheorem{prop}{\textbf{Proposition}}

\newcommand{\A}{\Bbb{A}}

\newcommand{\R}{\Bbb{R}}

\newcommand{\C}{\Bbb{C}}

\newcommand{\HH}{\Bbb{H}}

\newcommand{\Q}{\Bbb{Q}}

\newcommand{\F}{\Bbb{F}}

\newcommand{\G}{\Bbb{G}}

\newcommand{\Z}{\Bbb{Z}}

\newcommand{\GSp}{\text{GSp}}

\newcommand{\Sp}{\text{Sp}}

\newcommand{\U}{\text{U}}

\newcommand{\GO}{\text{GO}}

\newcommand{\GSO}{\text{GSO}}

\newcommand{\Frob}{\text{Frob}}

\newcommand{\End}{\text{End}}

\newcommand{\Hom}{\text{Hom}}

\newcommand{\tr}{\text{tr}}

\newcommand{\im}{\text{im}}

\newcommand{\Gal}{\text{Gal}}

\newcommand{\GL}{\text{GL}}

\newcommand{\SL}{\text{SL}}

\newcommand{\y}{\hspace{6pt}}

\title{{\bf{On Galois representations and Hilbert-Siegel modular forms}}}

\author{Claus M. Sorensen}

\begin{document}

\date{}

\maketitle

\begin{abstract}
This article is a spinoff of the book of Harris and Taylor [HT], in which they prove the local Langlands conjecture for $\GL(n)$, and its companion paper by Taylor and Yoshida [TY] on local-global compatibility.
We record some consequences in the case of genus two Hilbert-Siegel modular forms. In other words, we are concerned with cusp forms $\pi$ on $\GSp(4)$ over a totally real field. When $\pi$ is globally generic (that is, has a non-vanishing Fourier coefficient), and $\pi$ has a Steinberg component at some finite place, we associate a Galois representation compatible with the local Langlands correspondence for $\GSp(4)$ defined by Gan and Takeda in a recent preprint [GT]. Over $\Q$, for $\pi$ as above, this leads to a new realization of the Galois representations studied previously by Laumon, Taylor and Weissauer. We are hopeful that our approach should apply more generally, once the functorial lift to $\GL(4)$ is understood, and once the ``book project'' is completed. A concrete consequence of the above compatibility is the following special case of a conjecture stated in [SU]: If $\pi$ has nonzero vectors fixed by a non-special maximal compact subgroup at $v$, the corresponding monodromy operator at $v$ has rank at most one.\footnote{{\it{Keywords}}: Galois representations, Hilbert-Siegel modular forms, monodromy}
\footnote{{\it{2000 AMS Mathematics Classification}}: 11F33, 11F41, 11F70, 11F80}
\end{abstract}

\section{Introduction}

The interplay between modular forms and Galois representations has shown to be extremely fruitful in number theory. For example, the Hasse-Weil conjecture, saying that the $L$-function of an elliptic curve over $\Q$ has meromorphic continuation to $\C$, was proved in virtue of this reciprocity. Some of the most basic examples are Hilbert modular forms and Siegel modular forms, both very well-studied in the literature. In this paper, we study a mixture of these. First, we introduce the Siegel upper half-space (of complex dimension three):
$$
\mathcal{H}\overset{\text{df}}{=}\{\text{$Z=X+iY \in M_2(\C)$ symmetric, with $Y$ positive definite}\}.
$$
For a moment, we will view the symplectic similitude group $\GSp(4)$ as an affine group scheme over $\Z$, by choosing the non-degenerate alternating form to be
$$
\begin{pmatrix} 0 & I \\ -I & 0\end{pmatrix}.
$$
Later, beyond this introduction, we will switch to a skew-diagonal form. The similitude character is denoted by $c$ throughout. We then consider the subgroup $\GSp(4,\R)^+$ of elements with positive similitude. It acts on $\mathcal{H}$ in the standard way, by linear fractional transformations. More precisely, by the formula:
$$
\text{$gZ=(AZ+B)(CZ+D)^{-1}$, $\y$ $g=\begin{pmatrix} A & B \\ C & D\end{pmatrix}\in \GSp(4,\R)^+$, $\y$ $Z \in \mathcal{H}$.}
$$
The automorphy factor $CZ+D$ will be denoted by $j(g,Z)$ from now on. Next, we look at the $d$-fold product of this setup. In more detail, we will concentrate on a discrete subgroup $\Gamma$ inside $\GSp(4,\R)^{+d}$, and its diagonal action on $\mathcal{H}^d$. To exhibit examples of such $\Gamma$, we bring into play a totally real number field $F$, of degree $d$ over $\Q$, and label the real embeddings by $\sigma_i$. This ordering is used to identify $\GSp(4,\mathcal{O})^+$ with a discrete subgroup $\Gamma$. Here the plus signifies that we look at elements whose similitude is a totally positive unit in $\mathcal{O}$, the ring of integers in $F$. Finally, in order to define modular forms for $\Gamma$, we fix weights
$$
\text{$\underline{k}_i=(k_{i,1},k_{i,2})$, $\y$ $k_{i,1}\geq k_{i,2}\geq 3$, $\y$ $i=1,\ldots,d$.} 
$$
For each $i$, introduce the irreducible algebraic representation $\rho_{\underline{k}_i}$ of $\GL(2,\C)$,
$$
\text{Sym}^{k_{i,1}-k_{i,2}}(\C^2)\otimes {\det}^{k_{i,2}}.
$$
The underlying space of $\rho_{\underline{k}_i}$ is just a space of polynomials in two variables, homogeneous of a given degree. Following [Bai], we then define a Hilbert-Siegel modular form for $\Gamma$, with weights $\underline{k}_i$, to be a holomorphic vector-valued function
$$
f: \mathcal{H}^d\rightarrow \bigotimes_{i=1}^d \text{Sym}^{k_{i,1}-k_{i,2}}(\C^2)\otimes {\det}^{k_{i,2}}
$$
satisfying the following transformation property for every $\underline{Z}\in \mathcal{H}^d$ and $\gamma \in \Gamma$:
$$
f(\gamma\underline{Z})=\bigotimes_{i=1}^d\rho_{\underline{k}_i}(c(\gamma_i)^{-1} j(\gamma_i,Z_i))\cdot f(\underline{Z}).
$$
Such $f$ are automatically holomorphic at infinity, by the Koecher principle. The form is often assumed to be cuspidal, that is, it vanishes at infinity. These modular forms have rich arithmetic properties. To exploit them, it is useful to switch to an adelic setup and instead look at automorphic representations. For instance, this immediately gives rise to a workable Hecke theory. Thus, from now on in this paper, instead of $f$ we will focus on a cuspidal automorphic representation $\pi$ of $\GSp(4)$ over the totally real field $F$. The analogues of Hilbert-Siegel modular forms are those $\pi$ which are holomorphic discrete series at infinity. It is a basic fact that such $\pi$ do {\it{not}} admit Whittaker models. However, it is generally believed that the $L$-function of $\pi$ coincides with the $L$-function of a $\pi'$ which {\it{does}} admit a Whittaker model. Hence, for our purposes, there is no serious harm in assuming the existence of such a model. For a while though, let us not make this assumption, and explain in detail the expectations regarding Hilbert-Siegel modular forms and their associated Galois representations. 

\medskip

\noindent Let $\pi$ be a cuspidal automorphic representation of $\GSp(4)$ over some totally real field $F$, and denote by $S_{\pi}$ the set of finite places where $\pi$ is ramified. For each infinite place $v$, we assume that $\pi_v$ is an essentially discrete series representation, and that $\pi_v$ has central character $a \mapsto a^{-w}$. Here $w$ is an integer, independent of $v$. Under these assumptions, and a choice of an isomorphism $\iota:\bar{\Q}_{\ell}\rightarrow \C$, it is expected that there should be a semisimple continuous Galois representation
$$
\rho_{\pi,\iota}:\Gal(\bar{F}/F)\rightarrow \GL_4(\bar{\Q}_{\ell})
$$
with the following properties: $\rho_{\pi,\iota}$ is unramified at $v \notin S_{\pi}$ not dividing $\ell$, and 
$$
L_v(s-\frac{3}{2},\pi,\text{spin})=\det(1-\iota\rho_{\pi,\lambda}(\text{Frob}_{v})\cdot q_v^{-s})^{-1}.
$$
Here $\text{Frob}_{v}$ is the geometric Frobenius. More prudently, the above spin $L$-factor should actually lie in $L[q_v^{-s}]$ for some number field $L$ inside $\C$, and instead of $\iota$ one could focus on the finite place of $L$ it defines. In the {\it{rational}} case $F=\Q$, the existence of $\rho_{\pi,\iota}$ is now known, due to the work of many people (Chai-Faltings, Laumon, Shimura, Taylor and Weissauer). See [Lau] and [Wei] for the complete result. For arbitrary $F$, not much is known. Of course, when $\pi$ is CAP (cuspidal associated to parabolic), or a certain functorial lift (endoscopy, base change or automorphic induction), $\rho_{\pi,\iota}$ is known to exist by [BRo]. However, in most of these cases $\rho_{\pi,\iota}$ is reducible. In the opposite case, that is, when $\pi$ genuinely belongs to $\GSp(4)$, the representation $\rho_{\pi,\iota}$ should be irreducible. Obviously, this is the case we are interested in. Actually, we will aim higher and consider the ramified places $S_{\pi}$ too. The impetus for doing so, is the recent work of Gan and Takeda [GT], in which they prove the local Langlands conjecture for $\GSp(4)$. To an irreducible admissible representation $\pi_v$, they associate an $L$-parameter  
$$
\text{rec}_{\text{GT}}(\pi_v):W_{F_v}'=W_{F_v}\times \SL(2,\C) \rightarrow \GSp(4,\C).
$$
The notation $\text{rec}_{\text{GT}}$ is ours. This correspondence is natural in a number of respects. For example, it preserves the $L$- and $\epsilon$-factors defined by Shahidi in the generic case [Sha]. See section 2.2.2 below for a discussion of the complete list of desiderata. Now, the representation $\rho_{\pi,\iota}$ should satisfy local-global compatibility. That is, for any finite place $v$ (not dividing $\ell$), the restriction $\rho_{\pi,\iota}|_{W_{F_v}}$ should correspond to $\text{rec}_{\text{GT}}(\pi_v)$ through the usual dictionary [Tat]. As it stands, this is only morally true; one has to twist $\pi_v$. The precise folklore prediction is:

\medskip

\noindent {\bf{Conjecture.}} 
{\it{Let $\pi$ be a cuspidal automorphic representation of $\GSp(4)$ over some totally real field $F$. Assume there is a cuspidal automorphic representation of $\GL(4)$ over $F$, which is a weak lift of $\pi$. Moreover, we assume that
$$
\text{$\pi^{\circ}\overset{\text{df}}{=}\pi\otimes|c|^{\frac{w}{2}}$ is unitary, for some $w \in \Z$.}
$$
Finally, at each infinite place $v$, we assume that $\pi_v$ is an essentially discrete series representation with the same central and infinitesimal character as the finite-dimensional irreducible algebraic representation $V_{\mu(v)}$ of highest weight
$$
\text{$t=\begin{pmatrix}t_1 & & & \\ & t_2 & & \\ & & t_3 & \\ & & & t_4\end{pmatrix}\mapsto t_1^{\mu_1(v)}t_2^{\mu_2(v)}c(t)^{\delta(v)-w}$, $\y$ $\delta(v)\overset{\text{df}}{=}\frac{1}{2}(w-\mu_1(v)-\mu_2(v))$.}
$$
Here $\mu_1(v)\geq \mu_2(v) \geq 0$ are integers such that $\mu_1(v)+\mu_2(v)$ has the same parity as $w$.
In particular, the central character $\omega_{\pi_v}$ is of the form $a \mapsto a^{-w}$ at each infinite place $v$. Under these assumptions, for each choice of an isomorphism $\iota:\bar{\Q}_{\ell}\rightarrow \C$, there is a unique irreducible continuous representation
$$
\rho_{\pi,\iota}:\Gal(\bar{F}/F)\rightarrow \GSp_4(\bar{\Q}_{\ell})
$$
characterized by the following property: For each finite place $v\nmid \ell$ of $F$, we have
$$
\iota\text{WD}(\rho_{\pi,\iota}|_{W_{F_v}})^{F-ss}\simeq \text{rec}_{\text{GT}}(\pi_v\otimes |c|^{-\frac{3}{2}}).
$$
Moreover, $\pi^{\circ}$ is tempered everywhere. Consequently, $\rho_{\pi,\iota}$ is pure of weight ${\bf{w}}\overset{\text{df}}{=}w+3$. The representation $\rho_{\pi,\iota}$ has the following additional properties:

\begin{itemize}
\item $\rho_{\pi,\iota}^{\vee}\simeq \rho_{\pi,\iota}\otimes \chi^{-1}$ where the similitude $\chi=\omega_{\pi^{\circ}}\cdot\chi_{\text{cyc}}^{-{\bf{w}}}$ is totally odd.

\item The representation $\rho_{\pi,\iota}$ is potentially semistable at any finite place $v|\ell$. Moreover, $\rho_{\pi,\iota}$ is crystalline at a finite place $v|\ell$ when $\pi_v$ is unramified.

\item The Hodge-Tate weights are given by the following recipe: Fix an infinite place $v$, and use the same notation for the place above $\ell$ it defines via $\iota$. 
$$
\dim_{\bar{\Q}_{\ell}}\text{gr}^j(\rho_{\pi,\iota}\otimes_{F_v}B_{dR})^{\Gal(\bar{F}_v/F_v)}=0,
$$
unless $j$ belongs to the set
$$
\delta(v)+\{0, \mu_2(v)+1, \mu_1(v)+2, \mu_1(v)+\mu_2(v)+3\},
$$
in which case the above dimension is equal to one. 
\end{itemize}
}}

\medskip

\noindent The notation used will be explained carefully in the main body of the text below. Our main result is a proof of this conjecture in a substantial number of cases:

\medskip

\noindent {\bf{Theorem A.}} 
{\it{The above conjecture holds for globally generic $\pi$ such that, for some finite place $v$, the local component $\pi_v$ is an unramified twist of the Steinberg representation.}}

\medskip

\noindent The proof is an application of the monumental work of Harris-Taylor [HT], and its refinement by Taylor-Yoshida [TY]. Let us briefly sketch the simple strategy: First, since $\pi$ is globally generic, one can lift it to an automorphic representation $\Pi$ on $\GL(4)$ using theta series, by utilizing the close connection with $\GO(3,3)$. This is a well-known, though unpublished, result of Jacquet, Piatetski-Shapiro and Shalika. Other proofs exist in the literature. For example, see [AS] for an approach using the converse theorem.  We make use of Theorem 13.1 in [GT], saying that the lift $\pi \mapsto \Pi$ is strong. That is, compatible with the local Langlands correspondence everywhere. We note that, by the Steinberg assumption, $\Pi$ must be cuspidal. Next, we base change $\Pi$ to a CM extension $E$ over $F$, and twist it by a suitable character $\chi$ to make it conjugate self-dual. We can now apply [HT] and [TY] to the representation $\Pi_E(\chi)$, in order to get a Galois representation $\rho_{\Pi_E(\chi),\iota}$ over $E$. Since $E$ is arbitrary, a delicate patching argument shows how to descend this collection to $F$, after twisting by $\rho_{\check{\chi},\iota}$.

\medskip

\noindent  We note that Theorem A continues to hold if some $\pi_v$ is a generalized Steinberg representation of Klingen type (see section 2 below), or a supercuspidal not coming from $\GO(2,2)$. 
The point being that the local lift $\Pi_v$ on $\GL(4)$ should remain a discrete series. However, eventually the book project of the Paris 7 GRFA seminar should make any local assumption at $v$ superfluous. See Expected Theorem 2.4 in [Har]. Furthermore, the assumption that $\pi$ is globally generic is used exclusively to get a strong lift to $\GL(4)$. Our understanding is that the current state of the trace formula should at least give a {\it{weak}} lift more generally. See [Art] in conjunction with [Whi]. In this respect, there is a very interesting preprint of Weissauer [We2], in which he proves that if $\pi$ is a discrete series at infinity, it is weakly equivalent to a globally generic representation. However, apparently he needs to work over $\Q$. Perhaps ideas from [Lab] will be useful in treating $F$ of degree at least two. In any case, to get a {\it{strong}} lift, one would have to show that the $L$-packets defined in [GT] satisfy the expected character relations. This seems to be quite difficult. 

\medskip

\noindent Next, we will record a few corollaries of Theorem A. Part of our original motiviation for writing this paper, was to determine the rank of the monodromy operator for $\rho_{\pi,\iota}$ at a place $v$ where $\pi_v$ is Iwahori-spherical (that is, has nonzero vectors fixed by an Iwahori-subgroup). The question is completely answered by:

\medskip

\noindent {\bf{Theorem B.}} 
{\it{Let $\rho_{\pi,\iota}$ be the Galois representation attached to a globally generic cusp form $\pi$ as in Theorem A. Let $v \nmid \ell$ be a finite place of $F$ such that $\pi_v$ is Iwahori-spherical and ramified. Then $\rho_{\pi,\iota}|_{I_{F_v}}$ acts unipotently. Moreover,
\begin{itemize}
\item $\pi_v$ of Steinberg type $\Longleftrightarrow$ monodromy has rank $3$. 
\item $\pi_v$ has a unique $J_Q$-fixed line $\Longleftrightarrow$ monodromy has rank $2$.
\item $\pi_v$ para-spherical $\Longleftrightarrow$ monodromy has rank $1$.
\end{itemize}
}}
\medskip

\noindent Here $J_Q$ denotes the Klingen parahoric, and by $\pi_v$ being {\it{para-spherical}} we mean that it has nonzero vectors fixed by a non-special maximal compact subgroup. The first two consequences are part of the Conjecture on p. 11 in [GTi], while the latter is part of Conjecture 3.1.7 on p. 41 in [SU]. The holomorphic analogue of the last consequence would have applications to the Bloch-Kato conjecture for modular forms of square-free level. See [SU] and the authors thesis [Sor].

\medskip

\noindent The proof of Theorem B essentially follows from local-global compatibility, and the classification of the Iwahori-spherical representations of $\GSp(4)$. The dimensions of the parahoric fixed spaces for each class of representations were tabulated in [Sch]. We only need the generic ones. For those six types of representations, the data are reproduced in Table A in section 4.5 below.

\medskip

\noindent At the other extreme, in the supercuspidal case, we obtain the following: 

\medskip

\noindent {\bf{Theorem C.}} 
{\it{Let $\rho_{\pi,\iota}$ be the Galois representation attached to a globally generic cusp form $\pi$ as in Theorem A. Let $v \nmid \ell$ be a finite place of $F$ such that $\pi_v$ is supercuspidal, and not a lift from $\GO(2,2)$. Then $\rho_{\pi,\iota}|_{W_{F_v}}$ is irreducible. Furthermore, $\rho_{\pi,\iota}$ is trivial on some finite index subgroup of $I_{F_v}$, and we have:
$$
\frak{f}_{\text{Swan}}(\rho_{\pi,\iota}|_{I_{F_v}})=4 \cdot \text{depth}(\pi_v).
$$
}}

\noindent Here $\frak{f}_{\text{Swan}}$ denotes the Swan conductor, closely related to the more commonly used Artin conductor, and the depth of $\pi_v$ is defined in [MP]. The precise definitions are recalled below in section 4.5. The proof of Theorem C relies on two essential ingredients. One is a formula, due to Bushnell and Frolich [BF], relating the depth to the conductor in the case of supercuspidals on $\GL(n)$.
The second is a paper of Pan, showing that the local theta correspondence preserves depth [Pan]. In particular, if $\pi_v$ has depth zero, we deduce that $\rho_{\pi,\iota}$ is tamely ramified at $v$. We finish the paper with another criterion for tame ramification, due to Genestier and Tilouine [GTi] over $\Q$: Suppose $\pi_v^{J_Q,\chi}$ is nonzero, for some non-trivial character $\chi$ of $\F_v^*$ inflated to the units, then $\rho_{\pi,\iota}$ is tamely ramified:
$$
\rho_{\pi,\iota}|_{I_{F_v}}=1 \oplus 1 \oplus \chi \oplus \chi.
$$
Here $\chi$ is the character of $I_{F_v}$ obtained via local class field theory. Moreover, one can arrange for the two eigenspaces, for $1$ and $\chi$, to be totally isotropic. 

\medskip

\noindent From our construction of the compatible system $\rho_{\pi,\iota}$ in Theorem A, one deduces that it is {\it{motivic}} in the sense defined on p. 60 in [BRo]: There is a smooth projective variety $X_{/F}$, and an integer $n$, such that $\rho_{\pi,\iota}$ is a constituent of
$$
H^j(X \times_F \bar{F},\bar{\Q}_{\ell})(n),
$$
for all $\ell$, where ${\bf{w}}=j-2n$. By invoking the Weil restriction, it is enough to show the analogous result for $\rho_{\pi,\iota}|_{\Gamma_E}$ for some CM extension $E$ over $F$. For a detailed argument, we refer to the proof of Proposition 5.2.1 on p. 86 in [BRo]. Over $E$, the variety is a self-product of the universal abelian variety over a simple Shimura variety. See the bottom isomorphism on p. 98 in [HT]. As in the case of Hilbert modular forms [BRo], one would like to have actual motives over $E$ associated with $\pi$. Seemingly, one of the main obstacles in deriving this from [HT] is a multiplicity one issue for the unitary groups considered there: Is the positive integer $a$ in part (6) on p. 12 of [TY] in fact equal to one? Conjecturally, one should even have motives over $F$ attached to $\pi$. Even over $\Q$ this is not yet known. One problem is that the Hecke correspondences on a Siegel threefold do not extend to a given toroidal compactification. For a more thorough discussion of these matters, and a slightly different approach, see [H].

\medskip

\noindent Many thanks are due to D. Ramakrishnan for his suggestion that I should look at the Hilbert-Siegel case by passing to a CM extension, and for sharing his insights on many occasions. I am also grateful to M. Harris for useful correspondence regarding the patching argument in section 4.3. Finally, I am thankful to C. Skinner and A. Wiles for discussions relevant to this paper, and for their encouragement and support.

\section{Lifting to $\GL(4)$}

We will describe below how to transfer automorphic representations of $\GSp(4)$, of a certain type, to $\GL(4)$. Throughout, we work over a totally real base field $F$. Let us take $\pi$ to be a globally generic cuspidal automorphic representation of $\GSp(4)$, with central character $\omega_{\pi}$. We do $not$ assume it is unitary. For $v|\infty$, 
$$
\text{$\pi_v\simeq\pi_{\mu(v)}^W$, $\y$
$\mu_1(v)\geq\mu_2(v)\geq0$, $\y$ 
$\mu_1(v)+\mu_2(v)\equiv w$ (mod $2$).}
$$
The notation is explained more carefully below. Here we fix the integer $w$ such that $\omega_{\pi_v}$ takes 
nonzero $a \mapsto a^{-w}$ for all archimedean places $v$. In particular, 
$$
\text{$\pi^{\circ}\overset{\text{df}}{=}\pi \otimes |c|^{\frac{w}{2}}$ $is$ unitary.}
$$
Using theta series, one can then associate an automorphic representation $\Pi$ of $\GL(4)$ with the following properties: It has central character $\omega_{\pi}^2$, and satisfies\footnote{We normalize the isomorphism $W_{\R}^{\text{ab}}\simeq \R^*$ using the absolute value $|z|_{\C}\overset{\text{df}}{=}|z|^2$.} 

\begin{itemize}
\item $\Pi\otimes \omega_{\pi}^{-1}\simeq \Pi^{\vee}$.
\item For $v|\infty$, the $L$-parameter of $\Pi_v$  has the following restriction to $\C^*$,
$$
z \mapsto |z|^{-w}\cdot \begin{pmatrix} (z/\bar{z})^{\frac{\nu_1+\nu_2}{2}} & & & \\ & (z/\bar{z})^{\frac{\nu_1-\nu_2}{2}} & & \\ & & (z/\bar{z})^{-\frac{\nu_1-\nu_2}{2}} & \\ & & & (z/\bar{z})^{-\frac{\nu_1+\nu_2}{2}}\end{pmatrix},
$$
where $\nu_1=\mu_1+2$ and $\nu_2=\mu_2+1$ give the Harish-Chandra parameter of $\pi_v$. Here we suppress the dependence on $v$, and simply write $\mu_i=\mu_i(v)$.
\item $L(s,\Pi_v)=L(s,\pi_v,\text{spin})$, for finite $v$ such that $\pi_v$ is unramified.
\end{itemize}
In fact, Gan and Takeda have recently defined a local Langlands correspondence for $\GSp(4)$ such that the above lift is $strong$. That is, the $L$-parameters of $\pi_v$ and $\Pi_v$ coincide at $all$ places $v$. For later applications, we would like $\Pi$ to have a square-integrable component. Using table 2 on page 51 in [GT], we can ensure this by assuming the existence of a finite place $v_0$ where $\pi_{v_0}$ is of the form  
$$
\text{$\pi_{v_0}=
\begin{cases}
\text{St}_{\GSp(4)}(\chi)\\
\text{St}(\chi,\tau)
\end{cases}$ 
$\Longrightarrow$ $\y$
$\Pi_{v_0}=
\begin{cases}
\text{St}_{\GL(4)}(\chi)\\
\text{St}(\tau).
\end{cases}$}
$$
Also, $\Pi_{v_0}$ is supercuspidal if $\pi_{v_0}$ is a supercuspidal not coming from $\GO(2,2)$.
Most of the notation used here is self-explanatory, except maybe the symbol $\text{St}(\chi,\tau)$:
It denotes the generalized Steinberg representation, of Klingen type, associated to a supercuspidal $\tau$ on $\GL(2)$ and a non-trivial quadratic character $\chi$ such that $\tau\otimes \chi=\tau$. We refer to page 35 in [GT] for more details. As a bonus, the existence of such a place $v_0$ guarantees that $\Pi$ is $cuspidal$: Otherwise, $\pi$ is a theta lift from $\GO(2,2)$, but by [GT] the above $\pi_{v_0}$ do not participate here.  

\subsection{The archimedean case}

\subsubsection{Discrete series for $\GL(2,\R)$}

In this section, we briefly set up notation for the discrete series representations of $\GL(2,\R)$. Throughout we use Harish-Chandra parameters, as opposed to Blattner parameters. We will follow the notation of [Lau]. Thus, for each positive integer $n$, we let $\sigma_n$ be the unique (essentially) discrete series representation of $\GL(2,\R)$ which has the same central character and the same infinitesimal character as the finite-dimensional irreducible representation $\text{Sym}^{n-1}(\C^2)$.
$$
\sigma_n(\lambda)\overset{\text{df}}{=}\sigma_n \otimes |\det|^{\lambda},
$$
for each $\lambda \in \C$. More concretely, $\sigma_n$ is induced from the neutral component:
$$
\sigma_n=\text{Ind}_{\GL(2,\R)^+}^{\GL(2,\R)}(\sigma_n^+),
$$
where $\sigma_n^+$ is a certain representation of $\GL(2,\R)^+$ on the Hilbert space of 
$$
\text{$f:\mathcal{H}\rightarrow \C$ holomorphic, $\y$ $\|f\|^2=\int_{\mathcal{H}}|f(x+iy)|^2y^{n-1}dxdy<\infty$.}
$$
Here $\mathcal{H}$ denotes the upper half-plane in $\C$, and $\GL(2,\R)^+$ acts by the formula
$$
\sigma_n^+\begin{pmatrix}a & c \\ b & d\end{pmatrix}f(z)=(ad-bc)^n(cz+d)^{-n-1}f(\frac{az+b}{cz+d}).
$$
If instead of $z$ we use $\bar{z}$ in the automorphy factor, this also defines a representation $\sigma_n^-$ on the $anti$-holomorphic functions on $\mathcal{H}$. Then $\sigma_n$ can be thought of as the direct sum $\sigma_n^+\oplus \sigma_n^-$, where a non-trivial coset representative acts by reflection in the $y$-axis.
Its Jacquet-Langlands correspondent $\sigma_n^{\text{JL}}$ is simply $\text{Sym}^{n-1}(\C^2)$ viewed as a representation of the Hamilton quaternions $\HH^*$ embedded into $\GL(2,\C)$ in the standard fashion. 
Note that $\sigma_n$ is $not$ unitary, unless $n=1$. However, after a suitable twist it becomes unitary. Moreover,
$$
\Hom_{\GL(2,\R)^+}(\sigma_n^+(\frac{1-n}{2}),L_{\text{cusp}}^2(\R_+^*\Gamma\backslash \GL(2,\R)^+))
$$ 
can be identified with the space of weight $n+1$ elliptic cusp forms for the discrete subgroup $\Gamma$.
The weight $n+1$ is the Blattner parameter of $\sigma_n$, describing its minimal $K$-types. Up to isomorphism, $\sigma_n$ is easily seen to be invariant under twisting by the sign character of $\R^*$. Consequently, $\sigma_n$ is automorphically induced from a character on $\C^*$. More generally, its twist $\sigma_n(\lambda)$ is induced from
$$
z \mapsto |z|^{2\lambda-1}z^n=|z|^{n-1+2\lambda}(z/\bar{z})^{\frac{n}{2}}.
$$
We want to write down an $L$-parameter for $\sigma_n(\lambda)$. Thus, we let $W_{\R}$ denote the Weil group of $\R$, generated by $\C^*$ and an element $j$ such that $j^2=-1$ and $jz=\bar{z}j$ for all $z \in \C^*$.
To $\sigma_n(\lambda)$ is associated a conjugacy class of homomorphisms
$$
\phi_n(\lambda):W_{\R}\rightarrow \GL(2,\C)
$$
with semisimple images. By the above remarks, a concrete representative is:
$$
\text{$\phi_n(\lambda):z \mapsto |z|^{n-1+2\lambda}\cdot\begin{pmatrix}(z/\bar{z})^{\frac{n}{2}} &  \\  & (z/\bar{z})^{-\frac{n}{2}}\end{pmatrix}$, $\y$ 
$\phi_n(\lambda):j \mapsto \begin{pmatrix}  & 1 \\ (-1)^n & \end{pmatrix}$.}
$$
We note that the image of $\phi_n(\lambda)$ is bounded precisely when $\sigma_n(\lambda)$ is unitary. 

\subsubsection{Discrete series for $\GSp(4,\R)$}

We parametrize the discrete series representations of $\GSp(4,\R)$ in accordance with [Lau]. Throughout, we realize symplectic groups with respect to the form
$$
\text{$J=\begin{pmatrix} & S \\ -S & \end{pmatrix}$, $\y$ $S=\begin{pmatrix} & 1 \\ 1 & \end{pmatrix}$.}
$$
The similitude character is denoted by $c$. We take $B$ to be the Borel subgroup consisting of upper triangular matrices. The maximal torus $T$ is of the form  
$$
T=\{t=\begin{pmatrix}t_1 & & & \\ & t_2 & & \\ & & t_3 & \\ & & & t_4\end{pmatrix}: c(t)=t_1t_4=t_2t_3\}.
$$
We identify its group of rational characters $X^*(T)$ with the set of triples of integers $\mu=\mu_0\oplus(\mu_1,\mu_2)$, such that $\mu_1+\mu_2\equiv \mu_0$ (mod $2$), using the recipe:
$$
t^{\mu}\overset{\text{df}}{=}t_1^{\mu_1}t_2^{\mu_2}c(t)^{\frac{\mu_0-\mu_1-\mu_2}{2}}.
$$
Its restriction to the center $\G_m$ takes $a \mapsto a^{\mu_0}$. Inside $X^*(T)$ we have the cone of $B$-dominant weights $X^*(T)^+$ consisting of all tuples $\mu$ such that $\mu_1 \geq \mu_2 \geq 0$. By a fundamental result of Chevalley, the finite-dimensional irreducible algebraic representations are classified by their highest weights. For a $B$-dominant weight $\mu$ as above, we let $V_{\mu}$ be the corresponding algebraic representation of $\GSp(4)$. Its central character is given by $\mu_0$ as described above. To describe the infinitesimal character of $V_{\mu}$, we consider half the sum of $B$-positive roots:
$$
\delta\overset{\text{df}}{=}0 \oplus(2,1).
$$
The Harish-Chandra isomorphism identifies the center of the universal enveloping algebra, $Z(\frak{g})$, with the invariant symmetric algebra $\text{Sym}(\frak{t}_{\C})^W$. Under this isomorphism, the aforementioned infinitesimal character corresponds to 
$$
\nu\overset{\text{df}}{=}\mu+\delta=\mu_0 \oplus (\nu_1,\nu_2)=\mu_0 \oplus (\mu_1+2,\mu_2+1).
$$
Up to infinitesimal equivalence, there are precisely $two$ essentially discrete series representations of $\GSp(4,\R)$ with the same central character and the same infinitesimal character as $V_{\mu}$. Together, they form an $L$-packet,
$$
\{\pi_{\mu}^W,\pi_{\mu}^H\}.
$$
Here $\pi_{\mu}^W$ is the unique {\it{generic}} member, that is, it has a Whittaker model. The other member $\pi_{\mu}^H$ is {\it{holomorphic}}, and does not have a Whittaker model. They both have central character $a \mapsto a^{\mu_0}$ for real nonzero $a$. Another way to distinguish the two representations, is to look at their $(\frak{g},K)$-cohomology: The generic member $\pi_{\mu}^W$ has cohomology of Hodge type $(2,1)$ and $(1,2)$, whereas $\pi_{\mu}^H$ contributes cohomology of Hodge type $(3,0)$ and $(0,3)$. For example,
$$
H^{3,0}(\frak{g},K; \pi_{\mu}^H\otimes V_{\mu}^*)\simeq H^{0,3}(\frak{g},K; \pi_{\mu}^H\otimes V_{\mu}^*)\simeq \C,
$$
and similarly for $\pi_{\mu}^W$. We have consistently used the Harish-Chandra parameter $\nu$. Other authors prefer the Blattner parameter, because of its connection to the weights of Siegel modular forms. In our case, the relation is quite simple:
$$
\text{$\underline{k}=(k_1,k_2)$, $\y$ $k_1=\nu_1+1=\mu_1+3$, $\y$ $k_2=\nu_2+2=\mu_2+3$.}
$$
See Theorem 12.21 in [Kn] for example. A word of caution: It is really only fair to call $\underline{k}$ the Blattner parameter in the holomorphic case. In the generic case, it does not give the highest weight of the minimal $K$-type. The restriction of $\pi_{\mu}^H$ to the neutral component $\GSp(4,\R)^+$, consisting of elements with positive similitude, decomposes as a direct sum of a holomorphic part $\pi_{\mu}^{H+}$ and a dual anti-holomorphic part $\pi_{\mu}^{H-}$. If we assume $\mu_0=0$, so that these are unitary,
$$
\Hom_{\GSp(4,\R)^+}(\pi_{\mu}^{H+}, L_{\text{cusp}}^2(\R_+^*\Gamma\backslash \GSp(4,\R)^+))
$$
can be identified with cuspidal Siegel modular forms of weight $\underline{k}$ for the discrete subgroup 
$\Gamma$. The minimal $K$-type of $\pi_{\mu}^{H+}$, where $K$ is isomorphic to $\U(2)$, is 
$$
\text{Sym}^{k_1-k_2}(\C^2)\otimes \text{det}^{k_2}.
$$
This is the algebraic representation of $\GL(2,\C)$ with highest weight $t_1^{k_1}t_2^{k_2}$. Next, we wish to explain how $\pi_{\mu}^W$ and $\pi_{\mu}^H$ can be described explicitly as certain theta lifts: Up to equivalence, there are precisely two $4$-dimensional quadratic spaces over $\R$ of discriminant one. Namely, the anisotropic space $V_{4,0}$ and the split space $V_{2,2}$. The latter can be realized as $\text{M}(2,\R)$ equipped with the determinant. The former can be taken to be $\HH$ endowed with the reduced norm. In particular,
$$
\text{$\GSO(2,2)=(\GL(2,\R)\times \GL(2,\R))/\R^*$, $\y$
$\GSO(4,0)=(\HH^*\times \HH^*)/\R^*$}.
$$
Here $\R^*$ is embedded in the centers by taking $a$ to the element $(a,a^{-1})$. Thus irreducible representations of $\GSO(2,2)$ correspond to pairs of irreducible representations of $\GL(2,\R)$ having the same central character. Similarly for $\GSO(4,0)$. Now suppose $\sigma$ and $\sigma'$ are irreducible representations of $\GL(2,\R)$ with the same central character. We say that the representation $\sigma\otimes \sigma'$ is $regular$ if $\sigma\neq \sigma'$. The induced representation of the whole similitude group $\GO(2,2)$ then remains irreducible, and we denote it by $(\sigma\otimes\sigma')^+$ following the notation of [Rob]. On the other hand, in the $invariant$ case where $\sigma=\sigma'$, there are exactly two extensions of $\sigma\otimes \sigma$ to a representation of $\GO(2,2)$. By Theorem 6.8 in [Rob], precisely one of these extensions participates in the theta correspondence with $\GSp(4,\R)$. It is again denoted by $(\sigma\otimes\sigma)^+$. The analogous results hold in the anisotropic case. Now, the following two identities can be found several places in the literature. See [Moe] for example, or Proposition 4.3.1 in [HK]:
$$
\pi_{\mu}^W=\theta((\sigma_{\nu_1+\nu_2}(\frac{1}{2}(\mu_0-\nu_1-\nu_2+1))
\otimes\sigma_{\nu_1-\nu_2}(\frac{1}{2}(\mu_0-\nu_1+\nu_2+1)))^+),
$$
and
$$
\pi_{\mu}^H=\theta((\sigma_{\nu_1+\nu_2}^{\text{JL}}(\frac{1}{2}(\mu_0-\nu_1-\nu_2+1))
\otimes\sigma_{\nu_1-\nu_2}^{\text{JL}}(\frac{1}{2}(\mu_0-\nu_1+\nu_2+1)))^+).
$$
Here the equalities signify infinitesimal equivalence. Since we are dealing with theta correspondence for $similitude$ groups, it is unnecessary to specify an additive character. As a result, we can exhibit a parameter for the $\mu$-packet above:
$$
\phi_{\mu}: z \mapsto |z|^{\mu_0}\cdot \begin{pmatrix} (z/\bar{z})^{\frac{\nu_1+\nu_2}{2}} & & & \\ & (z/\bar{z})^{\frac{\nu_1-\nu_2}{2}} & & \\ & & (z/\bar{z})^{-\frac{\nu_1-\nu_2}{2}} & \\ & & & (z/\bar{z})^{-\frac{\nu_1+\nu_2}{2}}\end{pmatrix},
$$
and
$$
\phi_{\mu}: j \mapsto \begin{pmatrix} & & & 1 \\ & & 1 & \\ & (-1)^{\mu_0+1} & & \\ (-1)^{\mu_0+1} & & & \end{pmatrix}.
$$
Visibly, $\phi_{\mu}$ maps into the dual of the elliptic endoscopic group, consisting of
$$
\begin{pmatrix} a & & & b \\ & e & f & \\ & g & h & \\ c & & & d\end{pmatrix}=
\begin{pmatrix} 1 & & &  \\ &  & 1 & \\ &  &  & 1 \\  & 1 & & \end{pmatrix}
\begin{pmatrix} a & b & &  \\ c & d &  & \\ &  & e & f \\  & & g & h\end{pmatrix}
\begin{pmatrix} 1 & & &  \\ &  & 1 & \\ &  &  & 1 \\  & 1 & & \end{pmatrix}^{-1}
$$
such that $ad-bc$ equals $eh-fg$. Furthermore, the restriction of $\phi_{\mu}\otimes|\cdot|^{-3}$ to $\C^*$ is a direct sum of $distinct$ characters of the form $z \mapsto z^p\bar{z}^q$ for two {\it{integers}} $p$ and $q$ such that $p+q=\mu_0-3$. They come in pairs: If the type $(p,q)$ occurs, so does $(q,p)$. We say that $\pi_{\mu}^W$ and $\pi_{\mu}^H$ are regular algebraic (up to a twist).

\subsubsection{The Langlands classification for $\GL(4,\R)$}

The Langlands classification for $\GL(4,\R)$ describes all its irreducible admissible representations up to infinitesimal equivalence. The building blocks are the essentially discrete series $\sigma_n(\lambda)$, and the characters $\text{sgn}^n(\lambda)$ of the multiplicative group $\R^*$. The representations of $\GL(4,\R)$ are then constructed by parabolic induction. For example, start out with the representations $\sigma_n(\lambda)$ and $\sigma_{n'}(\lambda')$. We view their tensor product as a representation of the parabolic associated with the partition $(2,2)$, by making it trivial on the unipotent radical. Consider 
$$
\text{Ind}_{P_{(2,2)}}^{\GL(4,\R)}(\sigma_n(\lambda)\otimes\sigma_{n'}(\lambda')),
$$
where we use $normalized$ induction. Consequently, this is unitary when $\sigma_n(\lambda)$ and $\sigma_{n'}(\lambda')$ are both unitary. By interchanging their roles, we may assume that
$$
\text{Re}(\lambda) \geq \text{Re}(\lambda').
$$
In this case, the induced representation has a unique irreducible quotient. We denote it by 
$\sigma_n(\lambda)\boxplus\sigma_{n'}(\lambda')$ and call it the $isobaric$ sum. Its $L$-parameter is
$$
\phi=\phi_n(\lambda)\oplus \phi_{n'}(\lambda'):W_{\R}\rightarrow \GL(4,\C).
$$
We want to know when it descends to a parameter for $\GSp(4,\R)$. First, in the case where $n=n'$, it maps into the Levi subgroup of the Siegel parabolic:
$$
P=\{
\begin{pmatrix} A & \\ & c \cdot{^{\tau}A^{-1}}\end{pmatrix}
\begin{pmatrix} 1 & & x & y \\ & 1 & z & x \\ & & 1 & \\ & & & 1\end{pmatrix}
\}.
$$
Here $\tau$ is transposition with respect to the $skew$ diagonal. This means the packet for $\GSp(4,\R)$ should be obtained by induction from the Klingen parabolic,
$$
Q=
\{
\begin{pmatrix} t & & \\ & A & \\ & & t^{-1}\det(A)\end{pmatrix}
\begin{pmatrix} 1 & z & & \\ & 1 & & \\ & & 1 & -z \\ & & & 1\end{pmatrix}
\begin{pmatrix} 1 & & x & y \\ & 1 & & x \\ & & 1 & \\ & & & 1 \end{pmatrix}
\}.
$$
To be more precise, consider the following (unitarily) induced representation:
$$
\text{Ind}_Q^{\GSp(4,\R)}(\sigma_n(\lambda')\otimes |\cdot|^{2\lambda-2\lambda'}).
$$
When $\text{Re}(\lambda)$ is $strictly$ greater than $\text{Re}(\lambda')$ it has a unique irreducible quotient. At the other extreme, when $\lambda=\lambda'$ it decomposes into a direct sum of two $limits$ of discrete series. Secondly, in the case where $n \neq n'$ we try to conjugate $\phi$ into the dual of the elliptic endoscopic group. The determinant condition becomes:
\begin{equation}
\lambda-\lambda'=\frac{1}{2}(n'-n)\in \Z_+.
\end{equation}
If this is satisfied, the isobaric sum descends to a packet for $\GSp(4,\R)$, whose members can be constructed by theta correspondence as discussed above:
$$
\text{$\{\pi_{\mu}^W,\pi_{\mu}^H\}$, $\y$ $\mu_0=n-1+2\lambda$, $\y$ $\nu_1=\frac{1}{2}(n'+n)$, $\y$ $\nu_2=\frac{1}{2}(n'-n)$.}
$$

\subsection{The non-archimedean case}

\subsubsection{The local Langlands correspondence for $\GL(n)$}

We will quickly review the parametrization of the irreducible admissible representations of $\GL(n,F)$, up to isomorphism, where $F$ is a finite extension of $\Q_p$. We will suppress $F$ and denote this set by $\Pi(\GL(n))$. This parametrization was originally conjectured by Langlands, for any connected reductive group, and recently proved for $\GL(n)$ in [HT] and [Hen] by two different methods. We let $W_F$ be the Weil group of $F$. That is, the dense subgroup of the Galois group acting as integral powers of Frobenius on the residue field. It gets a topology by decreeing that the inertia group $I_F$ is open. From local class field theory,
$$
F^* \overset{\sim}{\rightarrow} W_F^{\text{ab}}.
$$ 
Here the isomorphism is normalized such that uniformizers correspond to lifts of the geometric Frobenius. It is used tacitly to identify characters of $F^*$ with characters of $W_F$. For arbitrary $n$, we consider the set $\Phi(\GL(n))$ consisting of conjugacy classes of continuous semisimple $n$-dimensional representations 
$$
\phi:W_F'=W_F \times \SL(2,\C) \rightarrow \GL(n,\C).
$$
The group $W_F'$ is sometimes called the Weil-Deligne group. The local Langlands correspondence is then a canonical collection of bijections $\text{rec}_n$, one for each $n$, 
$$
\text{rec}_n: \Pi(\GL(n)) \overset{1:1}{\rightarrow} \Phi(\GL(n))
$$  
associating an $L$-parameter $\phi_{\pi}$ to a representation $\pi$. It satisfies a number of natural properties, which in fact determine the collection uniquely. Namely, 

\begin{itemize}
\item The bijection $\text{rec}_1$ is given by local class field theory as above.
\item For any two $\pi \in \Pi(\GL(n))$ and $\sigma \in \Pi(\GL(r))$, we have equalities
$$
\begin{cases}
L(s,\pi\times \sigma)=L(s,\phi_{\pi}\otimes \phi_{\sigma}), \\
\epsilon(s,\pi\times \sigma,\psi)=\epsilon(s,\phi_{\pi}\otimes \phi_{\sigma},\psi).
\end{cases}
$$
\item The $L$-parameter of $\pi\otimes\chi\circ\det$ equals $\phi_{\pi}\otimes \chi$, for any character $\chi$.
\item For any $\pi$ as above, $\det(\phi_{\pi})$ corresponds to its central character $\omega_{\pi}$.
\item For any $\pi$ as above, $\text{rec}_n(\pi^{\vee})$ is the contragredient of $\text{rec}_n(\pi)$.
\end{itemize}
Here $\psi$ is a non-trivial character of $F$, used to define the $\epsilon$-factors. The collection 
$\text{rec}_n$ does not depend on it. The $L$ and $\epsilon$-factors on the left-hand side are those from  
[JPS]: They are first defined for generic representations, such as supercuspidals, and then one extends the definition to all representations using the Langlands classification. For the explicit formulas, see [Kud] or [Wed]. For the right-hand side, the definition of the $L$ and $\epsilon$-factors can be found in [Tat]:
The $L$-factors are given fairly explicitly, whereas the $\epsilon$-factors are defined very implicitly. One only has an abstract characterization due to Deligne and Langlands. For a good review of these definitions, we again refer to [Wed]. It is useful to instead consider the (Frobenius semisimple) Weil-Deligne representation of $W_F$ associated with a parameter $\phi$ as above. This is a pair $(r,N)$ consisting of a semisimple representation $r$ of $W_F$, and an operator $N$ satisfying the equation
$$
r(w)\circ N \circ r(w)^{-1}=|w|_F\cdot N
$$
for all $w \in W_F$. This $N$ is called the monodromy operator, and it is automatically nilpotent. The correspondence relies on the Jacobson-Morozov theorem:
$$
\text{$r(w)=\phi(w,\begin{pmatrix} |w|_F^{\frac{1}{2}} & \\ & |w|_F^{-\frac{1}{2}}\end{pmatrix})$, $\y$ $\exp(N)=\phi(1,\begin{pmatrix}1 & 1 \\ & 1\end{pmatrix})$.}
$$
The local Langlands correspondence $\text{rec}_n$ satisfies a number of additional natural properties, expected to hold more generally, of which we mention only a few:

\begin{itemize}
\item $\pi$ is supercuspidal $\Leftrightarrow$ $\phi_{\pi}$ is irreducible (and the monodromy is trivial). 
\item $\pi$ is essentially discrete series $\Leftrightarrow$ $\phi_{\pi}$ does $not$ map into a proper Levi.
\item $\pi$ is essentially tempered $\Leftrightarrow$ $\phi_{\pi}|_{W_F}$ has $bounded$ image in $\GL_n(\C)$.
\item $\pi$ is generic $\Leftrightarrow$ the adjoint $L$-factor $L(s, \text{Ad}\circ \phi_{\pi})$ has $no$ pole at $s=1$.
\end{itemize}
Here $\text{Ad}$ denotes the adjoint representation of $\GL_n(\C)$ on its Lie algebra $\frak{gl}_n(\C)$.

\subsubsection{The local Langlands correspondence for $\GSp(4)$}

For $\GSp(4)$, the presence of endoscopy makes the parametrization of $\Pi(\GSp(4))$ more complicated. 
It is partitioned into finite subsets $L_{\phi}$, called $L$-packets, each associated with a parameter $\phi$ as above mapping into the subgroup $\GSp(4,\C)$. We use the notation $\Phi(\GSp(4))$ for the set of $\GSp(4,\C)$-conjugacy classes of such $\phi$. The first attempt to define these $L$-packets, when $p$ is $odd$, is the paper [Vig]. The crucial case is when $\phi$ does not map into a proper Levi subgroup. In this case, Vigneras defined certain subsets $L_{\phi}$ by theta lifting from various forms of $\GO(4)$. However, she did not prove that these $L_{\phi}$ $exhaust$ all of $\Pi(\GSp(4))$. The work of Vigneras was later refined, so as to include the case $p=2$, in the paper [Ro2]. More recently, Gan and Takeda [GT] were able to prove the exhaustion, for all primes $p$. To do that, they used work of Muic-Savin, Kudla-Rallis, and Henniart. The main theorem of [GT] gives a finite-to-one surjection
$$
L:\Pi(\GSp(4)) \twoheadrightarrow \Phi(\GSp(4)),
$$
attaching an $L$-parameter $\phi_{\pi}$ to a representation $\pi$, and having the properties:

\begin{itemize}
\item $\pi$ is essentially discrete series $\Leftrightarrow$ $\phi_{\pi}$ does $not$ map into a proper Levi. 
\item For any generic $or$ non-supercuspidal $\pi \in \Pi(\GSp(4))$, and $\sigma \in \Pi(\GL(r))$, 
$$
\begin{cases}
\gamma(s,\pi\times \sigma,\psi)=\gamma(s,\phi_{\pi}\otimes \phi_{\sigma},\psi),\\
L(s,\pi\times \sigma)=L(s,\phi_{\pi}\otimes \phi_{\sigma}), \\
\epsilon(s,\pi\times \sigma,\psi)=\epsilon(s,\phi_{\pi}\otimes \phi_{\sigma},\psi).
\end{cases}
$$
\item The $L$-parameter of $\pi\otimes\chi\circ c$ equals $\phi_{\pi}\otimes \chi$, for any character $\chi$.
\item For any $\pi$ as above, $c(\phi_{\pi})$ corresponds to its central character $\omega_{\pi}$.
\end{itemize}
In the generic case, the invariants occurring on the left-hand side of the second condition are those from [Sha]. The definition can be extended to non-generic non-supercuspidals, using the Langlands classification. See page 13 in [GT]. For non-generic supercuspidals, $L$ satisfies an additional technical identity,
which we will not state here. It expresses a certain Plancherel measure as a product of four $\gamma$-factors. One has to include this last property to ensure the $uniqueness$ of $L$, as long as a satisfying theory of $\gamma$-factors is absent in this setup. For completeness, let us mention a few extra properties of the map $L$: For a given parameter $\phi$, the elements of the fiber $L_{\phi}$ correspond to characters of the group
$$
A_{\phi}=\pi_0(Z_{\GSp(4,\C)}(\im\phi)/\C^*)=\begin{cases}
\Z/2\Z, \\
0.
\end{cases}
$$
Moreover, an $L$-packet $L_{\phi}$ contains a generic member exactly when $L(s,\text{Ad}\circ \phi)$ has no pole at $s=1$. If in addition $\phi|_{W_F}$ has bounded image, the members of $L_{\phi}$ are all essentially tempered, and the generic member is unique. It is indexed by the trivial character of $A_{\phi}$. Next, we wish to at least give some idea of how the reciprocity map $L$ is constructed in [GT]: The key is to make use of theta liftings from various orthogonal similitude groups. In analogy with the archimedean case, there are two $4$-dimensional quadratic spaces over $F$ of discriminant one. We abuse notation slightly, and continue to denote the anisotropic space by $V_{4,0}$ and the split space by $V_{2,2}$. They can be realized as $D$ equipped with the reduced norm, where $D$ is a possibly split quaternion algebra over $F$. Again,
$$
\text{$\GSO(2,2)=(\GL(2,F)\times \GL(2,F))/F^*$, $\y$
$\GSO(4,0)=(D^*\times D^*)/F^*$},
$$  
as previously, where $D$ is here the $division$ quaternion algebra. Furthermore, one looks at the $6$-dimensional quadratic space $D \oplus V_{1,1}$. When $D$ is split, this is simply $V_{3,3}$. There is then a natural isomorphism, as given on page 9 in [GT],
$$
\GSO(3,3)=(\GL(4,F)\times F^*)/\{(a\cdot I,a^{-2}): a \in F^*\}.
$$
Now, start with an irreducible representation $\pi$ of $\GSp(4,F)$. By Theorem 5.3 in [GT], which relies on the work [KR] of Kudla and Rallis on the conservation conjecture, it follows that there are two possible $mutually$ $exclusive$ scenarios: 

\begin{enumerate}
\item $\pi$ participates in the theta correspondence with $\GSO(4,0)$,
\item $\pi$ participates in the theta correspondence with $\GSO(3,3)$.
\end{enumerate}
In the $first$ case, one has two essentially discrete series representations $\sigma$ and $\sigma'$ of $\GL(2,F)$ having the same central character, such that $\pi$ is the theta lift
$$
\pi=\theta((\sigma^{\text{JL}} \otimes \sigma'^{\text{JL}})^+)=\theta((\sigma'^{\text{JL}} \otimes \sigma^{\text{JL}})^+).
$$
By the local Langlands correspondence for $\GL(2)$, we have associated parameters $\phi_{\sigma}$ and 
$\phi_{\sigma'}$ with equal determinants. We then conjugate their sum $\phi_{\sigma}\oplus\phi_{\sigma'}$ into the dual of the elliptic endoscopic group as in the archimedean case,
$$
\phi_{\pi}=\phi_{\sigma}\oplus\phi_{\sigma'}:W_F' \rightarrow \GL(2,\C)\times_{\C^*} \GL(2,\C) \subset \GSp(4,\C).
$$
In the $second$ case, we write $\theta(\pi)$ as a tensor product $\Pi \otimes \omega_{\pi}$ for an irreducible representation $\Pi$ of $\GL(4,F)$. The local Langlands correspondence for $\GL(4)$ yields a parameter $\phi_{\Pi}$. We need to know that it maps into $\GSp(4,\C)$ after conjugation. When $\pi$ is a $discrete$ series, this follows from a result of Muic and Savin [MS], stated as Theorem 5.4 in [GT]: Indeed, the exterior square $L$-factor
$$
L(s,\Pi,\wedge^2\otimes \omega_{\pi}^{-1})
$$ 
has a pole at $s=0$. When $\pi$ is $not$ a discrete series, Gan and Takeda compute $\theta(\pi)$ explicitly,
using standard techniques developed by Kudla. For a summary, we refer to Table 2 on page 51 in their paper [GT]. It follows by {\it{inspection}} that $\phi_{\Pi}$ can be conjugated into a Levi subgroup of $\GSp(4,\C)$. Their computation works even for $p=2$, and hence completes the {\it{exercise}} of Waldspurger [Wal]. Finally, in Proposition 11.1 of [GT], it is shown by a global argument that the above construction is consistent with that of Vigneras and Roberts. 

\subsection{The globally generic case}

In the global situation, functoriality predicts that one should be able to transfer automorphic representations from $\GSp(4)$ to $\GL(4)$. It is widely believed that this should eventually follow by using trace formula techniques. See [Art] for a discussion on this approach. In the $globally$ $generic$ case, it has been known for some time that one can obtain (weak) lifts using theta series. This was first announced by Jacquet, Piatetski-Shapiro and Shalika, but to the best of our knowledge they never wrote it up.
However, many of the details are to be found in [Sou]. Moreover, there is an alternative proof in [AS] relying on the converse theorem. In this section, we wish to quote a recent refinement of the above transfer, due to Gan and Takeda [GT]. First, for completeness, let us recall the notion of being globally generic: Consider the upper-triangular Borel subgroup  
$$
B=\{\begin{pmatrix} s & & & \\ & t & & \\ & & ct^{-1} & \\ & & & cs^{-1}\end{pmatrix}
\begin{pmatrix} 1 & u & & \\ & 1 & & \\ & & 1 & -u \\ & & & 1 \end{pmatrix}
\begin{pmatrix} 1 & & x & y \\ & 1 & z & x \\ & & 1 & \\ & & & 1\end{pmatrix}
\}.
$$
We let $N$ denote its unipotent radical. Now let $F$ be a number field, and pick a non-trivial character $\psi$ of $\A_F$ trivial on $F$. By looking at $\psi(u+z)$ we can view it as an automorphic character of $N$.
An automorphic representation $\pi$ of $\GSp(4)$ over $F$ is then said to be globally generic if the Whittaker functional
$$
\text{$f \mapsto \int_{N(F)\backslash N(\A_F)}f(n)\psi^{-1}(n)dn$, $\y$ $f \in \pi$,}
$$
is not identically zero. This notion does not depend on $\psi$. As a consequence of the theta series approach we are about to discuss, Soudry proved the strong multiplicity one property in [Sou] for globally generic cusp forms $\pi$ on $\GSp(4)$. As mentioned above, the exterior square $\wedge^2$ defines an isogeny between the groups $\GL(4)$ and $\GSO(3,3)$. Thus, we need the global theta correspondence for {\it{similitude}} groups. For this, we refer to section 5 of [HK], noting that the normalization there differs slightly from [GT]. We will quickly review the main features of the definition: The Weil representation $\omega_{\psi}$ extends naturally to
$$
R=\{(g,h) \in \GSp(4)\times \GO(3,3): c(g)\cdot c(h)=1\}.
$$
Then, a Schwartz-Bruhat function $\varphi$ defines a theta kernel $\theta_{\varphi}$ on $R$ by the usual formula. The theta series lifting of a form $f$ in $\pi$ is hence given by the integral
$$
\theta_{\varphi}(f)(h)=\int_{\Sp(4,F)\backslash \Sp(4,\A_F)}\theta_{\varphi}(gg_h,h)f(gg_h)dg,
$$
where $g_h$ is any element with inverse similitude $c(h)$. The space spanned by all such theta series $\theta_{\varphi}(f)$ constitute an automorphic representation of $\GO(3,3)$, which we will denote by $\theta(\pi)$. It is independent of $\psi$. By Proposition 1.2 in [Sou], it is nonzero precisely when $\pi$ is globally generic. In fact, one can express the Whittaker functional for $\theta(\pi)$ in terms of that for $\pi$ given above. In particular, $\theta(\pi)$ is always generic, even though it may not be cuspidal. 
From now on, we will only view $\theta(\pi)$ as a representation of the subgroup $\GSO(3,3)$. As such, it remains irreducible. See Lemma 3.1 in [GT] for example. In turn, via the identification $\wedge^2$ we view
$\theta(\pi)$ as a representation $\Pi \otimes \omega_{\pi}$. We then have:

\begin{thm}
The global theta lifting $\pi \mapsto \theta(\pi)=\Pi \otimes \omega_{\pi}$ defines an injection from the set of globally generic cuspidal automorphic representations $\pi$ of $\GSp(4)$ to the set of generic automorphic representations $\Pi$ of $\GL(4)$ with central character $\omega_{\Pi}=\omega_{\pi}^2$. Moreover, this lifting has the following properties:

\begin{itemize}
\item $\Pi\simeq \Pi^{\vee}\otimes \omega_{\pi}$. 
\item It is a {\bf{strong}} lift, that is, $\phi_{\Pi_v}=\phi_{\pi_v}$ for all places $v$.
\item The image of the lifting consists {\bf{precisely}} of those $\Pi$ satisfying:
\begin{enumerate}
\item $\Pi$ is cuspidal and $L^{S}(s,\Pi,\wedge^2\otimes \omega_{\pi}^{-1})$ has a pole at $s=1$, {\bf{or}}
\item $\Pi=\sigma \boxplus\sigma'$ for cuspidal $\sigma\neq \sigma'$ on $\GL(2)$ with central character $\omega_{\pi}$.
\end{enumerate}
In the latter case, $\pi$ is the theta lift of the cusp form $\sigma \otimes \sigma'$ on $\GSO(2,2)$.
\end{itemize}
\end{thm}

\noindent $Proof$. See section 13 in [GT]. $\square$

\medskip

\noindent The refinements, due to Gan and Takeda, are primarily: The characterization of the image, and the fact that the global lift is compatible with the local Langlands correspondence at all (finite) places. The result that the lift is $strong$, in this sense, essentially follows from the construction of the local reciprocity map. However, locally one has to check that $\phi_{\Pi}$ is equivalent to $\phi_{\sigma}\oplus \phi_{\sigma'}$, when $\pi$ is the theta lift of $any$ $\sigma\otimes \sigma'$ on $\GSO(2,2)$. This is the content of Corollary 12.13 in [GT]. It is a result of their explicit determination of the theta correspondence.   

\section{Base change to a CM extension}

In this section we will construct representations, for which the results from [HT] on Galois representations apply. For that purpose, we will fix an arbitrary CM extension $E/F$, and an arbitrary Hecke character $\chi$ of $E$ with the property:
$$
\chi|_{\A_F^*}=\omega_{\pi}^{-1}.
$$
Since $\omega_{\pi_v}$ is of the form $a \mapsto a^{-w}$ at each infinite place $v$, it follows that every such $\chi$ is automatically algebraic. For all but finitely many $E$, the global base change 
$\Pi_E$ is $cuspidal$. We twist it by $\chi$ to make it conjugate self-dual:
$$
\Pi_E(\chi)^{\theta}\simeq \Pi_E(\chi)^{\vee}.
$$
Furthermore, $\Pi_E(\chi)$ is regular algebraic of weight zero, and we can arrange for it to have at least one square integrable component, by imposing the condition that $E$ splits completely at $v_0$. Hence Theorem C in [HT] applies. In this section, we briefly review results of Arthur and Clozel on base change for $\GL(n)$, and discuss the compatibility with the local Langlands correspondence. 

\subsection{Local base change}

Now that the $p$-adic local Langlands correspondence is available for $\GL(n)$, due to the works of Harris-Taylor and Henniart, base change makes sense for an arbitrary finite extension of local fields $E/F$. 
Indeed, a representation $\Pi_E$ of $\GL(n,E)$ is the base change of a representation $\Pi$ of $\GL(n,F)$ precisely when
$$
\phi_{\Pi_E}=\phi_{\Pi}|_{W_E'}.
$$
However, eventually we will use results from [AC]. At the time this book was written, one had to resort to a harmonic analytic definition of base change which we will review below. Fortunately, the compatibility of the two definitions has been checked by other authors. We will give precise references later. The latter definition only works for a cyclic extension $E/F$. For simplicity, we take it to be quadratic, and let $\theta$ be the non-trivial element in its Galois group. By lemma 1.1 in [AC], the norm map on $GL(n,E)$, 
taking $\gamma\mapsto\gamma\gamma^{\theta}$, defines an injection    
$$
\mathcal{N}:\{\text{$\theta$-conjugacy classes in $\GL(n,E)$}\}\hookrightarrow\{\text{conjugacy classes in $\GL(n,F)$}\}.
$$
This is used to define transfer of orbital integrals: Two compactly supported smooth functions $f$ and $f_E$ are said to have matching orbital integrals when
$$
O_{\gamma}(f)=\begin{cases}
\text{$TO_{\delta\theta}(f_E)$, $\y$ $\gamma=\mathcal{N}\delta$},\\
\text{$0$, $\y$ $\gamma$ is $not$ a norm}.
\end{cases}
$$
For the definitions of the integrals involved here, we refer to page 15 in [AC]. It is the content of Proposition 3.1 in [AC] that any $f_E$ has a matching function $f$. The fundamental lemma in this case is Theorem 4.5 in [AC]. We can now state Shintani's definition of local base change, following Definition 6.1 in [AC]: Let $\Pi$ and $\Pi_E$ be irreducible admissible representations of $\GL(n,F)$ and $\GL(n,E)$ respectively, and assume that $\Pi_E^{\theta}$ is isomorphic to $\Pi_E$. Let $I_{\theta}$ be an intertwining operator between these, normalized such that $I_{\theta}^2$ is the identity. This determines $I_{\theta}$ up to a sign. We then say that $\Pi_E$ is a base change of $\Pi$ if and only if 
$$
\tr(\Pi_E(f_E)\circ I_{\theta})=c \cdot\tr \Pi(f)
$$
for all matching functions $f$ and $f_E$ as above. The non-zero constant $c$ depends only on the choice of measures, and of $I_{\theta}$. By Theorem 6.2 in [AC] local base change makes sense for tempered representations. Using the Langlands classification, the lift then extends to all representations. For this, see the discussion on page 59 in [AC]. Since Shintani's definition is employed in [AC], we will need: 

\begin{thm}
Shintani's harmonic analytic definition of the local cyclic base change lifting is compatible with the local Langlands correspondence for $\GL(n)$.
\end{thm}

\noindent $Proof$. In the non-archimedean case, this is part 5 of Lemma VII.2.6 on page 237 in [HT]. 
The archimedean case was settled, in general, by Clozel [Clo]. $\square$

\medskip

\noindent As an example in the archimedean case, let us base change $\sigma_n(\lambda)\boxplus \sigma_{n'}(\lambda')$ to $\GL(4,\C)$. For simplicity, we will stick to the case of interest in this paper where it descends to a discrete series $L$-packet for $\GSp(4,\R)$. That is, we assume (1).
$$
-w\overset{\text{df}}{=}n-1+2\lambda=n'-1+2\lambda'
$$
Following [Kna] we let $[z]$ denote $\frac{z}{|z|}$. Then the unitarily induced representation
$$
\text{Ind}_B^{\GL(4,\C)}([\cdot]^n \otimes [\cdot]^{-n} \otimes [\cdot]^{n'}\otimes [\cdot]^{-n'})\otimes |\det|^{-w}
$$
has a unique irreducible quotient. This is the base change we are looking for. The Langlands correspondence for $\GL(n,\C)$, which is much simpler than the real case, was first studied by Zelobenko and Naimark. A good reference for their results is the expository paper [Kna]. In the non-archimedean case, 
$$
\text{$\text{St}_{\GL(4)}(\chi)_E=\text{St}_{\GL(4)}(\chi_E)$, $\y$ $\chi_E=\chi\circ N_{E/F}$.}
$$
However, the generalized Steinberg representation $\text{St}(\tau)$ may $not$ base change to a discrete series. Indeed, if $\omega_{E/F}$ denotes the associated quadratic character,
$$
\text{$\text{St}(\tau)_E=\begin{cases}
\text{$\text{St}(\tau_E)$}\\
\text{$\text{St}_{\GL(2)}(\psi)\boxplus \text{St}_{\GL(2)}(\psi^{\theta})$}
\end{cases}$when $\y$ 
$\begin{cases}
\tau\neq\tau \otimes \omega_{E/F}\\
\tau=\tau \otimes \omega_{E/F}.
\end{cases}$}
$$
Here $\psi\neq \psi^{\theta}$ is a certain character of $E^*$, with automorphic induction $\tau$.

\subsection{Global base change}

We now let $E/F$ denote an arbitrary CM extension of the totally real field $F$. Thus, the extension $E/F$ is quadratic, and $E$ is totally imaginary. We let $\theta$ be the non-trivial element in the Galois group.
Let $\Pi$ and $\Pi_E$ be automorphic representations of $\GL(n,\A_F)$ and $\GL(n,\A_E)$ respectively. We will assume $\Pi_E$ is invariant under $\theta$. Then, we say that $\Pi_E$ is a $strong$ base change lift of $\Pi$ if
$$
\Pi_{E,w}=\Pi_{v,E_w} 
$$
for $all$ places $w|v$. Here, the right-hand side is the local base change of $\Pi_v$ to $\GL(n,E_w)$. 
When $v$ is $split$ in $E$, this lift is naturally identified with $\Pi_v$. By comparing trace formulas for $\GL(n)$, no stabilization required, Arthur and Clozel proved that such lifts always exist. More precisely, we have the following:

\begin{thm}
There is a unique {\bf{strong}} base change lift $\Pi \mapsto \Pi_E$ between isobaric automorphic representations on $\GL(n,\A_F)$ and $\GL(n,\A_E)$, satisfying: 

\begin{itemize}
\item $\omega_{\Pi_E}=\omega_{\Pi}\circ N_{E/F}$. 
\item The image of the lifting consists {\bf{precisely}} of the $\theta$-invariant $\Pi_E$.
\item If $\Pi$ is cuspidal, $\Pi_E$ is cuspidal if and only if $\Pi\neq \Pi\otimes \omega_{E/F}$.
\end{itemize}
\end{thm}

\noindent $Proof$. This is essentially Theorem 4.2 combined with Theorem 5.1 in [AC]. Arthur and Clozel makes the assumption that $\Pi$ is $induced$ from $cuspidal$. This is now superfluous; the residual spectrum of $\GL(n)$ is understood by [MW]. $\square$

\medskip

\noindent Note that, if $\Pi$ is cuspidal, $\Pi_E$ is cuspidal for all but finitely many CM extensions $E$.  
Indeed, the self-twist condition is satisfied if the discriminant of $E$ does not divide the conductor of $\Pi$. The theory of base change goes hand-in-hand with $automorphic$ $induction$, which is a strong lift 
from isobaric automorphic representations of $\GL(n,\A_E)$ to those of $\GL(2n,\A_F)$ compatible with the $\text{rec}_n$,
$$
\text{$\pi \mapsto I_E^F(\pi)$, $\y$ $\phi_{I_E^F(\pi)_v}=\text{Ind}_{E_w}^{F_v}(\phi_{\pi_w})$,}
$$
for all $w|v$. Again, this is due to Arthur and Clozel in much greater generality. See Theorem 6.2 in [AC]. For the compatibility with $\text{rec}_n$ at the ramified places, we again refer to Lemma VII.2.6 in [HT].
In analogy with the above, we have: 

\begin{thm}
There is a {\bf{strong}} automorphic induction lift $\pi \mapsto I_E^F(\pi)$ between isobaric automorphic representations on $\GL(n,\A_E)$ and $\GL(2n,\A_F)$, satisfying: 

\begin{itemize}
\item $\omega_{I_E^F(\pi)}=\omega_{\pi}|_{\A_F^*}$. 
\item The image consists {\bf{precisely}} of the $\Pi$ such that $\Pi=\Pi\otimes \omega_{E/F}$.
\item If $\pi$ is cuspidal, $I_E^F(\pi)$ is cuspidal if and only if $\pi \neq \pi^{\theta}$.
\end{itemize}
\end{thm}

\noindent $Proof$. This is due to Arthur and Clozel. See section 6 of chapter 3 in [AC]. $\square$

\medskip

\noindent For the sake of completeness, let us mention a few links between base change and automorphic induction. For any two cuspidal $\pi$ and $\Pi$ as above, we have:
$$
\text{$I_E^F(\Pi_E)=\Pi\boxplus (\Pi\otimes\omega_{E/F})$, $\y$ $I_E^F(\pi)_E=\pi \boxplus \pi^{\theta}$.}
$$

\subsection{Conjugate self-dual twists}

Let us now take any CM extension $E/F$, and consider $\Pi_E$, where $\Pi$ is the theta series lifting of our original globally generic $\pi$ on $\GSp(4)$. We will assume $\pi_{v_0}$ is of (twisted) Steinberg type at some finite place $v_0$ of $F$. Hence $\Pi_{E,w_0}$ is Steinberg at all places $w_0$ of $E$ dividing $v_0$, and this ensures $\Pi_E$  is cuspidal. In [HT], one associates Galois representations to certain {\it{conjugate self-dual}} representations. $\Pi_E$ itself may not satisfy this condition, when $\omega_{\pi,E}\neq 1$, but certain twists $do$:
$$
\text{$\Pi_E(\chi)^{\theta}\simeq \Pi_E(\chi)^{\vee}$ $\Leftrightarrow$ $\chi|_{\A_F^*}=\omega_{\pi}^{-1}\omega_{E/F}^{n}$,}
$$
for $n=0$ or $n=1$. Such Hecke characters $\chi$ of $E$ exist: Indeed, by Frobenius reciprocity, any 
Hecke character of $F$ has infinitely many extensions to $E$; they are precisely the constituents of the induced representation of the compact idele class group $C_E^1$. By modifying this argument slightly, one can even control the ramification of the extensions if need be. Now, recall that for every place $v|\infty$,
$$
\omega_{\pi_v}(a)=a^{-w},
$$
for all $a \in \R^*$. Therefore, to retain algebraicity, take $n=0$ above. In this case, it follows that $all$ the extensions $\chi$ are automatically algebraic. That is,
$$
\text{$\chi_w(z)=z^{a}\bar{z}^{b}$, $\y$ $a=a(w)\in \Z$, $\y$ $b=b(w) \in \Z$, $\y$ $a+b=w$, }
$$
for each infinite place $w$ of $E$, not to be confused with the weight! For such characters $\chi$, the twist $\Pi_E(\chi)$ remains regular algebraic, and the weight is $zero$.

\section{Galois representations}

\subsection{Galois representations over CM extensions}

Let $E$ be an arbitrary CM extension of the totally real field $F$. One of the ultimate goals of the $book$ project [Har], is to attach an $\ell$-adic Galois representation $\rho_{\Pi,\iota}$ to a regular algebraic conjugate self-dual cuspidal automorphic representation $\Pi$ of $\GL(n,\A_E)$, and a choice of an isomorphism $\iota: \bar{\Q}_{\ell}\rightarrow \C$. See $expected$ Theorem 2.4 in [Har] for a more precise formulation. In the case where $\Pi$ has a square-integrable component at some finite place, pioneering work on this problem was done by Clozel [Cl] and Kottwitz [Kot], relating $\rho_{\Pi,\iota}|_{W_{E_w}}$ to the unramified component $\Pi_w$ at most places $w$. Their work was later extended to all places $w \nmid \ell$ in [HT], by Harris and Taylor, in the course of proving the local Langlands conjecture. However, in [HT] the $monodromy$ operator is ignored. This issue has been taken care of by Taylor and Yoshida in [TY], resulting in:

\begin{thm}
Let $E$ be a CM extension of a totally real field $F$, and let $\Pi$ be a cuspidal automorphic representation of $\GL(n,\A_E)$ satisfying the conditions:

\begin{itemize}
\item $\Pi_{\infty}$ is regular algebraic, $H^{\bullet}(\frak{g},K; \Pi_{\infty}\otimes \mathcal{V}^*)\neq 0$.  
\item $\Pi$ is conjugate self-dual, $\Pi^{\vee}\simeq \Pi^{\theta}$.
\item $\Pi_{w_0}$ is (essentially) square integrable for some finite place $w_0$.
\end{itemize}

\noindent Fix an isomorphism $\iota: \bar{\Q}_{\ell}\rightarrow \C$. Then there is a continuous  representation 
$$
\rho_{\Pi,\iota}: \Gal(\bar{F}/E)\rightarrow \GL(n,\bar{\Q}_{\ell})
$$
such that for every finite place $w$ of $E$, not dividing $\ell$, we have the following:
$$
\iota\text{WD}(\rho_{\Pi,\iota}|_{W_{E_w}})^{F-ss}\simeq \text{rec}_n(\Pi_w\otimes |\det|^{\frac{1-n}{2}}).
$$
\end{thm}

\noindent $Proof$. This is Theorem 1.2 in [TY], which is a refinement of Theorem VII.1.9 in [HT]. Indeed, the former result identifies the monodromy operator: By Corollary VII.1.11 in [HT] it is known that $\Pi_w$ is $tempered$ for all finite places $w$. Therefore, by parts (3) and (4) of Lemma 1.4 in [TY], it suffices to show that $\iota\text{WD}(\rho_{\Pi,\iota}|_{W_{E_w}})$ is $pure$. That is, up to a shift, the weight filtration coincides with the monodromy filtration. This is proved in [TY] by a careful study of the Rapoport-Zink weight spectral sequence, the main new ingredient being the vanishing outside the middle-degree in Proposition 4.4 in [TY]. $\square$

\medskip

\noindent A word about the notation used in the previous Theorem: First, $\mathcal{V}$ denotes an irreducible algebraic representation over $\C$ of the group $R_{E/\Q}\GL(n)$, which we will denote by $\mathcal{G}$. Then $\frak{g}$ denotes the Lie algebra of $\mathcal{G}(\R)$, and $K$ is a maximal compact subgroup times $Z_{\mathcal{G}}(\R)$. The symbol $WD(\rho)$ stands for the Weil-Deligne representation corresponding to an $\ell$-adic representation $\rho$ of $W_{E_w}$, where $w \nmid \ell$. This pair $(r,N)$ is obtained by fixing a lift $\Frob_w$ of the geometric Frobenius, and a continuous surjective homomorphism $t_{\ell}:I_{E_w}\rightarrow \Z_{\ell}$, and then writing
$$
\rho(\Frob_w^m\sigma)=r(\Frob_w^m\sigma)\text{exp}(t_{\ell}(\sigma)N)
$$
for $\sigma \in I_{E_w}$ and integers $m$. Then $r$ is a representation of $W_{E_w}$ having an $open$ kernel, and $N$ is a nilpotent operator satisfying the formula mentioned above:
$$
N \circ r(\Frob_w)=q_w\cdot r(\Frob_w) \circ N.
$$
Here $q_w$ is the order of the residue field of $E_w$. The isomorphism class of $(r,N)$ is independent of the choices made. Finally, the superscript $F-ss$ signifies Frobenius {\it{semisimplification}}. That is, leave $N$ unchanged, but semisimplify $r$. The representations $\rho_{\Pi,\iota}$ above satisfy a number of additional nice properties:

\begin{itemize}
\item $\Pi$ is square integrable at some finite $w_0\nmid\ell$ $\Longrightarrow$ $\rho_{\Pi,\iota}$ is irreducible.
\item Let $w\nmid \ell$ be a finite place of $E$, and let $\alpha$ be an eigenvalue of $\rho_{\Pi,\iota}(\sigma)$ for some $\sigma \in W_{E_w}$. Then $\alpha$ belongs to $\bar{\Q}$, and for every embedding $\bar{\Q}\hookrightarrow \C$,
$$
|\alpha| \in q_w^{\frac{\Z}{2}}.
$$
\item Let $w\nmid \ell$ be a finite place, with $\Pi_w$ {\it{unramified}}, and let $\alpha$ be an eigenvalue of $\rho_{\Pi,\iota}(\Frob_w)$. Then $\alpha$ belongs to $\bar{\Q}$, and for every embedding $\bar{\Q}\hookrightarrow \C$,
$$
|\alpha|=q_w^{\frac{n-1}{2}}.
$$
\item The representation $\rho_{\Pi,\iota}$ is potentially semistable at any finite place $w|\ell$. Moreover, $\rho_{\Pi,\iota}$ is crystalline at a finite place $w|\ell$ when $\Pi_w$ is unramified.
\end{itemize}
The first property is clear since $\text{rec}_n(\Pi_{w_0}\otimes |\det|^{\frac{1-n}{2}})$ is indecomposable.
This was observed in Corollary 1.3 in [TY]. It is expected to continue to hold if $w_0|\ell$ (and even without the square-integrability condition, admitting the book project). The second and the third property are parts 1 and 2 of Theorem VII.1.9 in [HT]. The former is a special case of Lemma I.5.7 in [HT], which apparently follows from the Rapoport-Zink weight spectral sequence in conjunction with de Jong's theory of alterations. The latter follows from Deligne's work on the Weil conjectures. The last property comes down to the comparison theorems of $p$-adic Hodge theory. To clarify these comments, we will briefly sketch how  
$\rho_{\Pi,\iota}$ is realized geometrically in [HT]: One starts off with an $n^2$-dimensional central division algebra $B$ over $E$, equipped with a positive involution $\ast$ such that $\ast|_{E}=\theta$. It is assumed to satisfy a list of properties, which are irrelevant for our informal discussion. For a fixed $\beta \in B$ such that $\beta\beta^*=1$, look at the unitary similitude group $G$ defined as follows: For a commutative $\Q$-algebra $R$,
$$
G(R)=\{\text{$x \in (B^{\text{op}}\otimes_{\Q}R)^*$: $x^*\beta x=c(x)\beta$, with $c(x)\in R^*$}\}.
$$
The element $\beta$ is chosen such that, at infinity, the derived group takes the form
$$
G^{\text{der}}(\R)=U(n-1,1)\times U(n)^{[F:\Q]-1}.
$$
The group $G$ has an associated Shimura variety of PEL type. That is, for each sufficiently small compact open subgroup $K$ inside the finite adeles $G(\A_f)$, there is a smooth proper variety $X_K$ over $E$ classifying isogeny classes of polarized abelian schemes $A$ of dimension $[F:\Q]n^2$, endowed with a certain homomorphism from $B$ into $\End(A)_{\Q}$ and a so-called level-structure relative to $K$. On $X_K$ one defines a $\bar{\Q}_{\ell}$-sheaf $\mathcal{L}_{\xi}$ by fixing an algebraic representation $\xi$ of $G$ defined over $\bar{\Q}_{\ell}$. Then, consider the following direct limit over subgroups $K$, endowed with natural commuting actions of $G(\A_f)$ and of the Galois group of $E$:
$$
H^m(X,\mathcal{L}_{\xi})\overset{\text{df}}{=}\underrightarrow{\lim}_K H_{\text{et}}^m(X_K \times_E\bar{F},\mathcal{L}_{\xi})=\bigoplus_{\pi_f}\pi_f \otimes R_{\xi}^m(\pi_f).
$$
Here $\pi_f$ runs over the irreducible admissible representations of $G(\A_f)$, and $R_{\xi}^m(\pi_f)$ is a finite-dimensional continuous Galois representation of $E$. To construct $\rho_{\Pi,\iota}$, we first descend $\Pi$ to an automorphic representation $\tilde{\Pi}$ of $B^{\text{op},*}$ via the Jacquet-Langlands correspondence. Using results of Clozel and Labesse, one shows that $\psi \times \tilde{\Pi}$ is a base change from $G$, for some algebraic Hecke character $\psi$ of $E$. In this way, we end up with an automorphic representation $\pi$ of $G$, and 
$$
R_{\xi}^{n-1}(\pi_f^{\vee})^{\text{ss}}\simeq \rho_{\Pi,\iota}^a \otimes \rho_{\psi,\iota}
$$
for some positive integer $a$. For details see p. 228 in [HT], and p. 12 in [TY].

\subsection{Hodge-Tate weights}

We will now describe the Hodge-Tate weights, which are certain numerical invariants of the restriction of $\rho_{\Pi,\iota}$ to the Galois group of $E_w$ for each place $w|\ell$. We briefly recall their definition: Let $\mathcal{K}$ be a finite extension of $\Q_{\ell}$. Following Fontaine, we introduce the field of $\ell$-adic periods $B_{dR}$. This is a $\mathcal{K}$-algebra, it comes equipped with a discrete valuation $v_{dR}$, and has residue field $\C_{\ell}$. The topology on $B_{dR}$ is coarser than the one coming form the valuation. The Galois group of $\mathcal{K}$ acts continuously on $B_{dR}$. If $t \in B_{dR}$ is a uniformizer, then we have
$$
g\cdot t=\chi_{\text{cyc}}(g)t,
$$
where $\chi_{\text{cyc}}$ is the $\ell$-adic cyclotomic character. The valuation defines a filtration:
$$
\text{$\text{Fil}^j(B_{dR})\overset{\text{df}}{=}t^jB_{dR}^+$, $\y$ $\text{gr}^j(B_{DR})\overset{\text{df}}{=}\text{Fil}^j(B_{dR})/\text{Fil}^{j+1}(B_{dR})\simeq \C_{\ell}(j)$.}
$$
If $V$ is a finite-dimensional continuous $\bar{\Q}_{\ell}$-representation of $\Gal(\bar{\Q}_{\ell}/\mathcal{K})$, we let
$$
D_{dR}(V)\overset{\text{df}}{=}(V \otimes_{\Q_{\ell}}B_{dR})^{\Gal(\bar{\Q}_{\ell}/\mathcal{K})}.
$$
This is a module over $\bar{\Q}_{\ell}\otimes_{\Q_{\ell}}\mathcal{K}$, inheriting a filtration from $B_{dR}$.
We say that $V$ is $de$ $Rham$ if this module is free of rank $\dim_{\bar{\Q}_{\ell}}(V)$. In this case, for each embedding $\tau:\mathcal{K}\rightarrow \bar{\Q}_{\ell}$ we introduce a multiset of integers $\text{HT}_{\tau}(V)$. It contains $\dim_{\bar{\Q}_{\ell}}(V)$ elements, and $j$ occurs with multiplicity equal to the dimension of 
$$
\text{gr}^j(V \otimes_{\tau,\mathcal{K}}B_{dR})^{\Gal(\bar{\Q}_{\ell}/\mathcal{K})}=\text{gr}^jD_{dR}(V)\otimes_{\bar{\Q}_{\ell}\otimes_{\Q_{\ell}}\mathcal{K}, 1 \otimes \tau}\bar{\Q}_{\ell}.
$$
over $\bar{\Q}_{\ell}$. If this is nonzero, $j$ is called a $Hodge$-$Tate$ weight for $V$ relative to the embedding $\tau$. Now we specialize the discussion, and take $V$ to be the restriction of $\rho_{\Pi,\iota}$ as above. We will quote a result from [HT], as stated in [Har], relating the Hodge-Tate weights to the highest weights of the algebraic representation $\mathcal{V}^*$. Recall, this is the irreducible algebraic representation of $\mathcal{G}(\C)$ such that the tensor product $\Pi_{\infty}\otimes \mathcal{V}^*$ has cohomology. By the definition of the $\Q$-group $\mathcal{G}$,
$$
\text{$\mathcal{G}(\R)=\prod_{\sigma \in \Sigma}\GL_n(E \otimes_{F,\sigma}\R)$, $\y$ 
$\mathcal{G}(\C)=\prod_{\sigma \in \Sigma}\GL_n(E \otimes_{F,\sigma}\C)$,}
$$
where $\Sigma$ is the set of embeddings $\sigma:F \rightarrow \R$. For each such $\sigma$, following [Har], we let $\{\tilde{\sigma},\tilde{\sigma}^c\}$ denote the two complex embeddings of $E$ extending it. We write
$$
\text{$\mathcal{V}^*=\otimes_{\sigma \in \Sigma}\mathcal{V}_{\sigma}^*$, $\y$ $\mathcal{V}_{\sigma}^*=\mathcal{V}_{\tilde{\sigma}}^*\otimes \mathcal{V}_{\tilde{\sigma}^c}^*$.}
$$
Here $\mathcal{V}_{\tilde{\sigma}}^*$ is naturally identified with an irreducible algebraic representation of the group $\GL_n(\C)$, and we consider its highest weight relative to the $lower$ triangular Borel.
This is the character of the diagonal torus corresponding to 
$$
\mu(\tilde{\sigma})=(\mu_1(\tilde{\sigma})\leq \mu_2(\tilde{\sigma})\leq \cdots \leq \mu_n(\tilde{\sigma})).
$$
Similarly, we get a dominant $n$-tuple of integers $\mu(\tilde{\sigma}^c)$ for $\mathcal{V}_{\tilde{\sigma}^c}^*$. It is given by:
$$
\mu_i(\tilde{\sigma}^c)=-\mu_{n-i+1}(\tilde{\sigma}),
$$
by the polarization condition. The multisets $\text{HT}_{\tau}$ for $\rho_{\Pi,\iota}$ are determined by:

\begin{thm}
Fix an embedding $s:E \rightarrow \bar{\Q}_{\ell}$, and let $w$ denote the associated finite place of $E$ above $\ell$. Then the Hodge-Tate weights of $\rho_{\Pi,\iota}$ restricted to the Galois group $\Gal(\bar{\Q}_{\ell}/E_w)$ at $w$, where $E_w=s(E)^-$, are all of the form
$$
\text{$j=i-\mu_{n-i}(\iota(s)^c)$, $\y$ $i=0,\ldots,n-1$.}
$$
In particular, the Hodge-Tate weights all occur with multiplicity one. 
\end{thm}

\noindent $Proof$. This is part 4 of Theorem VII.1.9 on p. 227 in [HT], but with the normalization used in [Har] in (2.6) on p. 5: The shift from $\iota(s)$ to $\iota(s)^c$ reflects the fact that we work with the $dual$ of the $\Pi$ in [HT]. Note that the inequalities on p. 3 in [Har] should be reversed. $\square$

\subsection{Patching}

The next key step is to $descend$ the family $\rho_{\Pi_E(\chi),\iota}\otimes \chi^{-1}$ to the base field $F$. This is done by a patching argument, used in various guises by other authors. For example, see Proposition 4.3.1 in [BRo], or section 4.3 in [BRa]. Here we will use a variant of Proposition 1.1 in [Har], which in turn is based on the discussion on p. 230-231 in [HT]. The proof in [Har] is somewhat brief, and somewhat imprecise at the end, so we decided to include a more detailed proof below. Hopefully, this might serve as a convenient reference. In this section, we use $\Gamma_F$ as shorthand notation for the absolute Galois group $\Gal(\bar{F}/F)$. The setup is the following: We let $\mathcal{I}$ be a set of cyclic Galois extensions $E$, of a fixed number field $F$, of $prime$ degree $q_E$. For every $E \in \mathcal{I}$ we assume we are given an $n$-dimensional continuous semisimple $\ell$-adic Galois representation over $E$,
$$
\rho_E:\Gamma_E \rightarrow \GL_n(\bar{\Q}_{\ell}).
$$
Here $\ell$ is a fixed prime. The family of representations $\{\rho_E\}$ is assumed to satisfy:

\begin{enumerate}
\item[(a)] \underline{Galois invariance}: $\rho_E^{\sigma}\simeq \rho_E$, $\y$ $\forall \sigma \in \Gal(E/F)$,
\item[(b)] \underline{Compatibility}: $\rho_E|_{\Gamma_{EE'}}\simeq \rho_{E'}|_{\Gamma_{EE'}}$, 
\end{enumerate}
for all $E$ and $E'$ in $\mathcal{I}$. These conditions are certainly necessary for the $\rho_E$ to be of the form $\rho|_{\Gamma_E}$ for a representation $\rho$ of $\Gamma_F$. What we will show, is that in fact (a) and (b) are also sufficient conditions if $\mathcal{I}$ is large enough. That is,

\begin{df}
Following [Har], for a finite set $S$ of places of $F$, we say that $\mathcal{I}$ is $S$-general if and only if the following holds: For any finite place $v \notin S$, and any finite extension $M$ of $F$, there is an $E \in \mathcal{I}$ linearly disjoint from $M$ such that $v$ splits completely in $E$. In this case, there will be infinitely many such $E$.
\end{df}

\noindent Recall that since $E$ is Galois over $F$, it is linearly disjoint from $M$ precisely when $E \cap M=F$. Moreover, it is of prime degree, so this just means $E$ is not contained in $M$. Hence, $\mathcal{I}$ being $S$-general is equivalent to: For $v \notin S$, there are infinitely many $E \in \mathcal{I}$ in which $v$ splits. A slightly stronger condition is:

\begin{df}
We say that $\mathcal{I}$ is strongly $S$-general if and only if the following holds: For any finite set $\Sigma$ of places of $F$, disjoint from $S$, there is an $E \in \mathcal{I}$ in which every $v \in \Sigma$ splits completely.
\end{df}

\noindent To see that this is indeed {\it{stronger}}, we follow Remark 1.3 in [Har]: Fix a finite place $v \notin S$, and a finite extension $M$ of $F$. Clearly we may assume $M \neq F$ is Galois. Let $\{M_i\}$ be the subfields of $M$, Galois over $F$, with a simple Galois group. For each $i$ we then choose a place $v_i$ of $F$, not in $S$, which does $not$ split completely in $M_i$. We take $\Sigma$ to be $\{v,v_i\}$ in the above definition, and get an $E \in \mathcal{I}$ in which $v$ and every $v_i$ splits. If $E$ was contained in $M$, it would be one of the $M_i$, but this contradicts the choice of $v_i$. Thus, $E$ and $M$ are disjoint.

\medskip
\noindent $Example$. Let $\Sigma=\{p_i\}$ be a finite set of primes. As is well-known, for odd $p_i$,
$$
\text{$p_i$ splits in $\Q(\sqrt{d})$ $\Longleftrightarrow$ $p_i \nmid d$ and $\begin{pmatrix}\frac{d}{p_i}\end{pmatrix}=1$.}
$$
Here $d$ is any square-free integer. Moreover, $2$ splits when $d \equiv 1$ mod $8$. The set of $all$ integers $d$ satisfying the congruences $d \equiv 1$ mod $p_i$, for all $i$, form an arithmetic progression. By Dirichlet's Theorem, it contains infinitely many primes. Therefore, the following family of imaginary quadratic extensions
$$
\mathcal{I}=\{\text{$\Q(\sqrt{-p})$: almost all primes $p$}\}
$$
is strongly $\varnothing$-general. This gives rise to a similar family of CM extensions of any given totally real field $F$, by taking the set of all the composite fields $F \mathcal{I}$.
\medskip

\noindent The main result of this section, is a strengthening of Proposition 1.1 in [Har]:

\begin{lem}
Let $\mathcal{I}$ be an $S$-general set of extensions $E$ over $F$, of prime degree $q_E$, and let $\rho_E$ be a family of semisimple Galois representations satisfying the conditions (a) and (b) above. Then there is a continuous semisimple 
$$
\text{$\rho: \Gamma_F \rightarrow \GL_n(\bar{\Q}_{\ell})$, $\y$ $\rho|_{\Gamma_E}\simeq \rho_E$,}
$$
for {\bf{all}} $E \in \mathcal{I}$. This determines the representation $\rho$ uniquely up to isomorphism.  
\end{lem}

\noindent $Proof$. The proof below is strongly influenced by the proofs of Proposition 1.1 in [Har], and of Theorem VII.1.9 in [HT]. We simply include more details and clarifications. The proof is quite long and technical, so we divide it into several steps. Before we construct $\rho$, we start off with noting that it is necessarily {\it{unique}}: Indeed, for any place $v \notin S$, we find an $E \in \mathcal{I}$ in which $v$ splits. In particular, $E_w=F_v$ for all places $w$ of $E$ dividing $v$. Thus, all the restrictions $\rho|_{\Gamma_{F_v}}$ are uniquely determined. We conclude that $\rho$ is unique, by the Cebotarev Density Theorem. For the construction of $\rho$, we first establish some {\it{notation}} used throughout the proof: 
We fix an arbitrary $base$ point $E_0 \in \mathcal{I}$, and abbreviate
$$
\text{$\rho_0\overset{\text{df}}{=}\rho_{E_0}$, $\y$ $\Gamma_0\overset{\text{df}}{=}\Gamma_{E_0}$, $\y$ $G_0\overset{\text{df}}{=}\Gal(E_0/F)$, $\y$ $q_0\overset{\text{df}}{=}q_{E_0}$.}
$$
We let $H$ denote the Zariski closure of $\rho_0(\Gamma_0)$ inside $\GL_n(\bar{\Q}_{\ell})$, and consider its identity component $H^{\circ}$. Define $M$ to be the finite Galois extension of $E_0$ with
$$
\text{$\Gamma_M=\rho_0^{-1}(H^{\circ})$, $\y$ $\Gal(M/E_0)=\pi_0(H)$.}
$$
Let $T$ be the set of isomorphism classes of irreducible constituents of $\rho_0$, ignoring multiplicities.
By property (a), the group $G_0$ acts on $T$ from the right. We note that $\tau$ and $\tau^{\sigma}$ occur in $\rho_0$ with the same multiplicity. We want to describe the $G_0$-orbits on $T$. First, we have the set $P$ of fixed points $\tau=\tau^{\sigma}$ for all $\sigma$. The set of non-trivial orbits is denoted by $C$. Note that any $c \in C$ has prime cardinality $q_0$. For each such $c$, we pick a representative $\tau_c \in T$, and let $C_0$ be the set of all these representatives $\{\tau_c\}$. Each $\tau \in C_0$ obviously has a trivial stabilizer in $G_0$.

\medskip
\noindent {\bf{Step 1}}: {\it{The extensions of $\rho_0$ to $\Gamma_F$.}} 

\medskip
\noindent Firstly, a standard argument shows that each $\tau \in P$ has an extension $\tilde{\tau}$ to $\Gamma_F$. This uses the divisibility of $\bar{\Q}_{\ell}^*$, in order to find a suitable intertwining operator $\tau\simeq\tau^{\sigma}$. All the other extensions are then obtained from $\tilde{\tau}$ as unique twists:
$$
\text{$\tilde{\tau}\otimes \eta$, $\y$ $\eta \in \hat{G}_0$}.
$$
Here $\hat{G}_0$ is the group of characters of $G_0$. Secondly, for a $\tau \in c$, we introduce
$$
\tilde{\tau}\overset{\text{df}}{=}\text{Ind}_{\Gamma_0}^{\Gamma_F}(\tau).
$$
Since $\tau$ is $not$ Galois-invariant, this $\tilde{\tau}$ is irreducible. It depends only on the orbit $c$ containing $\tau$, and it is invariant under twisting by $\hat{G}_0$. It has restriction
$$
\tilde{\tau}|_{\Gamma_0}\simeq {\bigoplus}_{\sigma \in G_0}\tau^{\sigma}.
$$
If we let $m_{\tau}$ denote the multiplicity with which $\tau \in T$ occurs in $\rho_0$, we get that
$$
\{\bigoplus_{\tau \in C_0}m_{\tau}\cdot \tilde{\tau}\}\oplus\{\bigoplus_{\tau \in P}{\bigoplus}_{\eta \in \hat{G}_0}m_{\tilde{\tau},\eta}\cdot (\tilde{\tau}\otimes \eta)\}
$$
is an extension of $\rho_0$ to $\Gamma_F$ for all choices of non-negative $m_{\tilde{\tau},\eta}\in \Z$ such that 
$$
{\sum}_{\eta \in \hat{G}_0} m_{\tilde{\tau},\eta}=m_{\tau}
$$
for every fixed $\tau \in P$.

\medskip
\noindent {\bf{Step 2}}: {\it{$\rho_0(\Gamma_{NE_0})$ is dense in $H$, when $N$ is linearly disjoint from $M$ over $F$.}} 

\medskip
\noindent To see this, let us momentarily denote the Zariski closure of $\rho_0(\Gamma_{NE_0})$ by $H_N$.
$N$ is a finite extension, so $H_N$ has finite index in $H$. Consequently, we deduce that $H_N^{\circ}=H^{\circ}$. Now, $NE_0$ and $M$ are linearly disjoint over $E_0$, and therefore
$$
\Gamma_0=\Gamma_{NE_0}\cdot\Gamma_M \Longrightarrow \rho_0(\Gamma_0)\subset \rho_0(\Gamma_{NE_0})\cdot H^{\circ}\subset H_N.
$$
Taking the closure, we obtain that $H_N=H$.

\medskip
\noindent {\bf{Step 3}}: {\it{If $N$ is a finite extension of $F$, linearly disjoint from $M$ over $F$. Then:

\begin{enumerate}
\item[(1)] $\tau|_{\Gamma_{NE_0}}$ is irreducible, for all $\tau \in T$. 
\item[(2)] $\tilde{\tau}|_{\Gamma_N}$ is irreducible, for all $\tau \in P$.
\item[(3)] $\tau|_{\Gamma_{NE_0}}\simeq\tau'|_{\Gamma_{NE_0}} \Rightarrow \tau \simeq \tau'$, for all $\tau,\tau' \in T$.
\item[(4)] $\tilde{\tau}|_{\Gamma_N}\simeq (\tilde{\tau}'\otimes \eta)|_{\Gamma_N} \Rightarrow \text{$\tau\simeq\tau'$ and $\eta=1$}$, for all $\tau,\tau' \in P$ and $\eta \in \hat{G}_0$.
\end{enumerate}}} 

\medskip
\noindent Parts (1) and (3) follow immediately from Step 2, and obviously (1) implies (2).
Also, part (3) immediately implies that $\tau \simeq \tau'$ in (4). Suppose $\eta \in \hat{G}_0$ satisfies:
$$
\tilde{\tau}|_{\Gamma_N} \simeq \tilde{\tau}|_{\Gamma_N} \otimes \eta|_{\Gamma_N}
$$
for some $\tau \in P$. In other words, $\eta|_{\Gamma_N}$ occurs in $\text{End}_{\Gamma_{NE_0}}(\tilde{\tau}|_{\Gamma_N})$, which is trivial by part (1). So, $\eta$ is trivial on $\Gamma_N$ and on $\Gamma_0$. Hence, $\eta=1$ by disjointness.

\medskip
\noindent {\bf{Step 4}}: {\it{$\tilde{\tau}|_{\Gamma_N}$ is irreducible for all $\tau \in C_0$. That is, part (2) holds for all $\tau \in T$.}} 

\medskip
\noindent Since $N$ and $E_0$ are linearly disjoint over $F$, we see that $\Gamma_F=\Gamma_E\cdot \Gamma_0$. Hence,
$$
\tilde{\tau}|_{\Gamma_N}=\text{Ind}_{\Gamma_0}^{\Gamma_F}(\tau)|_{\Gamma_N}\simeq \text{Ind}_{\Gamma_{NE_0}}^{\Gamma_N}(\tau|_{\Gamma_{NE_0}}),
$$
by Mackey theory. Now, $\tau|_{\Gamma_{NE_0}}$ is irreducible and $not$ Galois-invariant.

\medskip
\noindent {\bf{Step 5}}: {\it{Suppose $E \in \mathcal{I}$ is linearly disjoint from $M$ over $F$. Then, for a unique choice of non-negative integers $m_{\tilde{\tau},\eta,E}$ with $\eta$-sum $m_{\tau}$, we have the formula:}}
$$
\rho_E\simeq \{\bigoplus_{\tau \in C_0}m_{\tau}\cdot \tilde{\tau}|_{\Gamma_E}\}\oplus\{\bigoplus_{\tau \in P}{\bigoplus}_{\eta \in \hat{G}_0} m_{\tilde{\tau},\eta,E}\cdot (\tilde{\tau}\otimes \eta)|_{\Gamma_E}\}.
$$ 
{\it{In particular, $\rho_0$ and $\rho_E$ have a \underline{common} extension to $\Gamma_F$.}}

\medskip
\noindent The uniqueness of the $m_{\tilde{\tau},\eta,E}$ follows directly from part (4) in Step 3. Recall,
$$
\rho_E|_{\Gamma_{EE_0}}\simeq \rho_0|_{\Gamma_{EE_0}}\simeq 
\{\bigoplus_{\tau \in C_0}m_{\tau}\cdot {\bigoplus}_{\sigma \in G_0}\tau|_{\Gamma_{EE_0}}^{\sigma}\}\oplus\{\bigoplus_{\tau \in P}m_{\tau}\cdot \tau|_{\Gamma_{EE_0}}\},
$$
by the compatibility condition (b). Here all the $\tau|_{\Gamma_{EE_0}}^{\sigma}$ are distinct by (3).
First, let us pick an arbitrary $\tau \in P$. As representations of $G_0$, viewed as the Galois group of $EE_0$ over $E$ by disjointness, we have
$$
\Hom_{\Gamma_{EE_0}}(\tilde{\tau}|_{\Gamma_{E}},\rho_E)\simeq {\bigoplus}_{\eta \in \hat{G}_0}\dim_{\bar{\Q}_{\ell}}\Hom_{\Gamma_E}((\tilde{\tau}\otimes \eta)|_{\Gamma_E},\rho_E)\cdot \eta.
$$
The $\bar{\Q}_{\ell}$-dimension of the left-hand side clearly equals $m_{\tau}$, and the right-hand side defines the partition $m_{\tilde{\tau},\eta,E}$ of $m_{\tau}$. Next, let us pick an arbitrary $\tau \in C_0$.
By the same argument, using that $\tilde{\tau}$ is invariant under twisting by $\hat{G}_0$, we get:
$$
\Hom_{\Gamma_{EE_0}}(\tilde{\tau}|_{\Gamma_{E}},\rho_E)\simeq \dim_{\bar{\Q}_{\ell}}\Hom_{\Gamma_E}(\tilde{\tau}|_{\Gamma_E},\rho_E)\cdot {\bigoplus}_{\eta \in \hat{G}_0}\eta.
$$
Now the left-hand side obviously has dimension $m_{\tau}q_0$. We deduce that $\tilde{\tau}|_{\Gamma_E}$ occurs in $\rho_E$ with multiplicity $m_{\tau}$. Counting dimensions, we obtain the desired decomposition of $\rho_E$. Note that we have not used the Galois invariance of $\rho_E$. In fact, it is a {\it{consequence}} of the above argument, assuming $E\cap M=F$.

\medskip
\noindent {\bf{Step 6}}: {\it{Fix an $E_1 \in \mathcal{I}$ disjoint from $M$ over $F$. Introduce the representation}}
$$
\rho\overset{\text{df}}{=}\{\bigoplus_{\tau \in C_0}m_{\tau}\cdot \tilde{\tau}\}\oplus\{\bigoplus_{\tau \in P}{\bigoplus}_{\eta \in \hat{G}_0} m_{\tilde{\tau},\eta,E_1}\cdot (\tilde{\tau}\otimes \eta)\}.
$$
{\it{Then $\rho|_{\Gamma_E}\simeq \rho_E$ for all extensions $E \in \mathcal{I}$ linearly disjoint from $ME_1$ over $F$.}}

\medskip
\noindent By definition, and Step 5, we have that $\rho|_{\Gamma_{E_1}}\simeq \rho_{E_1}$. Take $E \in \mathcal{I}$ to be any extension, disjoint from $ME_1$ over $F$. We compare the decomposition of $\rho|_{\Gamma_E}$,
$$
\rho|_{\Gamma_E}=\{\bigoplus_{\tau \in C_0}m_{\tau}\cdot \tilde{\tau}|_{\Gamma_E}\}\oplus\{\bigoplus_{\tau \in P}{\bigoplus}_{\eta \in \hat{G}_0} m_{\tilde{\tau},\eta,E_1}\cdot (\tilde{\tau}\otimes \eta)|_{\Gamma_E}\},
$$
to the decomposition of $\rho_E$ in Step 5. We need to show the multiplicities match:
$$
\text{$m_{\tilde{\tau},\eta,E}=m_{\tilde{\tau},\eta,E_1}$, $\y$ $\forall \tau \in P$, $\y$ $\forall \eta \in \hat{G}_0$.}
$$
By property (b), for the pair $\{E,E_1\}$, we know that $\rho|_{\Gamma_E}$ and $\rho_E$ become isomorphic after restriction to $\Gamma_{EE_1}$. Once we prove $EE_1$ is linearly disjoint from $M$ over $F$, we are done by (2) and (4). The disjointness folllows immediately:
$$
EE_1 \otimes_F M \simeq E \otimes_F E_1 \otimes_F M\simeq E \otimes_F ME_1\simeq EE_1M.
$$
\medskip
\noindent {\bf{Step 7}}: {\it{$\rho|_{\Gamma_E}\simeq \rho_E$ for {\bf{all}} $E \in \mathcal{I}$.}}

\medskip
\noindent By Step 6, we may assume $E \in \mathcal{I}$ is {\it{contained}} in $ME_1$. Now take an auxiliary 
extension $\mathcal{E}\in \mathcal{I}$ linearly disjoint from $ME_1$ over $F$. Consequently, using (b),
$$
\rho|_{\Gamma_{\mathcal{E}}}\simeq \rho_{\mathcal{E}} \Rightarrow \rho|_{\Gamma_{E\mathcal{E}}}\simeq 
\rho_{\mathcal{E}}|_{\Gamma_{E\mathcal{E}}}\simeq \rho_E|_{\Gamma_{E\mathcal{E}}}.
$$
Thus, $\rho|_{\Gamma_E}$ agrees with $\rho_E$ when restricted to $\Gamma_{E\mathcal{E}}$. It suffices to show that the union of these subgroups $\Gamma_{E\mathcal{E}}$, as $\mathcal{E}$ varies, is dense in $\Gamma_E$. Again, we invoke the Cebotarev Density Theorem. Indeed, let $w$ be a place of $E$, lying above $v \notin S$. It is then enough to find an $\mathcal{E}\in \mathcal{I}$, as above, such that $w$ splits completely in $E\mathcal{E}$. Then $\Gamma_{E_w}$ is contained in $\Gamma_{E\mathcal{E}}$. We $know$, by the $S$-generality of $\mathcal{I}$, that we can find an $\mathcal{E}\in \mathcal{I}$, not contained in $ME_1$, in which $v$ splits completely. This $\mathcal{E}$ works: This follows from elementary splitting theory, as $E$ and $\mathcal{E}$ are disjoint.

\medskip

\noindent This finishes the proof of the patching lemma. $\square$

\medskip

\noindent $Remark$. From the proof above, we infer the following concrete description of the patch-up representation $\rho$. First fix $any$ $E_0 \in \mathcal{I}$, and let $P$ be the set of Galois-invariant constituents $\tau$ of $\rho_{E_0}$. For each such $\tau$, we fix an extension $\tilde{\tau}$ to $F$ once and for all. Furthermore, let $C_0$ be a set of representatives for the non-trivial Galois orbits of constituents of $\rho_{E_0}$. Then $\rho$ is of the following form
$$
\rho\simeq\{\bigoplus_{\tau \in C_0}m_{\tau}\cdot \text{Ind}_{\Gamma_0}^{\Gamma_F}(\tau)\}\oplus\{\bigoplus_{\tau \in P}{\bigoplus}_{\eta \in \Gal(E_0/F)^{\wedge}} m_{\tilde{\tau},\eta}\cdot (\tilde{\tau}\otimes \eta)\}.
$$
Here the $m_{\tilde{\tau},\eta}$ are $some$ non-negative integers with $\eta$-sum $m_{\tau}$, the multiplicity of $\tau$ in $\rho_{E_0}$. This fairly explicit description may be useful in deriving properties of $\rho$ from those of $\rho_{E_0}$.

\medskip

\noindent For future reference, we finish this section with a few remarks on the generalization of 
the patching lemma to $solvable$ extensions. Thus, $\mathcal{I}$ now denotes a collection of solvable Galois 
extensions $E$ over $F$, and we assume we are given Galois representations $\rho_E$, as above, satisfying (a) and (b). For any $L$ over $F$,   
$$
\mathcal{I}_L\overset{\text{df}}{=}\{E \in \mathcal{I}: L \subset E\}.
$$
Loosely speaking, we say that $\mathcal{I}$ is $S$-general if it is $S$-general in prime layers:

\begin{df}
Following [Har], for a finite set $S$ of places of $F$, we say that $\mathcal{I}$ is $S$-general if and only if the following holds: For every $L$ such that $\mathcal{I}_L\neq \varnothing$, 
$$
\{\text{prime degree extensions $K/L$, with $\mathcal{I}_K \neq \varnothing$}\}
$$
is $S(L)$-general in the previous sense. $S(L)$ denotes the places of $L$ above $S$.
\end{df}

\noindent From now on, we will make the additional hypothesis that all the extensions $E \in \mathcal{I}$ have uniformly {\it{bounded heights}}. That is, there is an integer $H_\mathcal{I}$
such that every index $[E:F]$ has at most $H_\mathcal{I}$ prime divisors (not necessarily distinct).

\begin{lem}
Assume the collection $\mathcal{I}$ has uniformly bounded heights. Then $\mathcal{I}$ is $S$-general if and only if the following condition holds for every $L$ with $\mathcal{I}_L\neq \varnothing$: Given a finite place $w \notin S(L)$ and a finite extension $M$ over $L$, there is an extension $E \in \mathcal{I}_L$ linearly disjoint from $M$ over $L$, in which $w$ splits completely.
\end{lem}

\noindent $Proof$. The {\it{if}} part follows immediately by unraveling the definitions. The {\it{only if}} part is proved by induction on the maximal height of the collection $\mathcal{I}_L$ over $L$, the height {\it{one}} case being the definition. Suppose $\mathcal{I}_L$ has maximal height $H$, and assume the lemma holds for smaller heights. Let $w$ and $M$ be as above. By $S$-generality, there is a prime degree extension $K$ over $L$ with $\mathcal{I}_K \neq \varnothing$, disjoint from $M$ over $L$, in which $w$ splits. 
Fix a place $\tilde{w}$ of $K$ above $w$. Now, $\mathcal{I}_K$ clearly has maximal height {\it{less}} than $H$. By the induction hypothesis there is an $E \in \mathcal{I}_K$,  disjoint from $MK$ over $K$, in which 
$\tilde{w}$ splits. This $E$ works. $\square$

\medskip

\noindent Under the above assumptions on $\mathcal{I}$, a given place $w \notin S(L)$ splits completely in infinitely many $E \in \mathcal{I}_L$, unless $L$ belongs to $\mathcal{I}$. One has a stronger notion:

\begin{df}
We say that $\mathcal{I}$ is strongly $S$-general if and only if the following holds: For any $L$ such that $\mathcal{I}_L\neq \varnothing$, and any finite set $\Sigma$ of places of $L$ disjoint from $S(L)$, there is an $E \in \mathcal{I}_L$ in which every $v \in \Sigma$ splits completely.
\end{df}

\noindent As in the prime degree case, treated above, one shows that this is indeed a stronger condition. Our next goal is to prove the following generalization of the patching lemma to certain collections of solvable extensions:

\begin{thm}
Let $\mathcal{I}$ be an $S$-general collection of solvable Galois extensions $E$ over $F$, with uniformly bounded heights, and let $\rho_E$ be a family of $n$-dimensional continuous semisimple $\ell$-adic Galois representations satisfying the conditions (a) and (b) above. Then there is a continuous semisimple representation 
$$
\text{$\rho: \Gamma_F \rightarrow \GL_n(\bar{\Q}_{\ell})$, $\y$ $\rho|_{\Gamma_E}\simeq \rho_E$,}
$$
for {\bf{all}} $E \in \mathcal{I}$. This determines the representation $\rho$ uniquely up to isomorphism.  
\end{thm}

\noindent $Proof$. Uniqueness is proved by paraphrasing the argument in the prime degree situation. The existence of $\rho$ is proved by induction on the maximal height of $\mathcal{I}$ over $F$, the height one case being the previous patching lemma. Suppose $\mathcal{I}$ has maximal height $H$, and assume the Theorem holds for smaller heights. Take an arbitrary prime degree extension $K$ over $F$, with $\mathcal{I}_K \neq \varnothing$. Clearly $\mathcal{I}_K$ is an $S(K)$-general set of solvable Galois extensions of $K$, of maximal height strictly smaller than $H$. Moreover, the subfamily $\{\rho_E\}_{E \in \mathcal{I}_K}$ obviously satisfies (a) and (b). By induction, we find a continuous semisimple $\ell$-adic representation 
$$
\text{$\rho_K: \Gamma_K \rightarrow \GL_n(\bar{\Q}_{\ell})$, $\y$ $\rho_K|_{\Gamma_E}\simeq \rho_E$,}
$$
for all $E \in \mathcal{I}_K$. We then wish to apply the prime degree patching lemma to the family $\{\rho_K\}$, as $K$ varies over extensions as above. By definition, such $K$ do form an $S$-general collection over $F$. It remains to show that $\{\rho_K\}$ satisfies (a) and (b). To check property (a), take any $\sigma \in \Gamma_F$, and note that $\rho_K^{\sigma}$ agrees with $\rho_K$ after restriction to $\Gamma_E$ for an arbitrary extension $E \in \mathcal{I}_K$. The union of these $\Gamma_E$ is dense in $\Gamma_K$ by the Cebotarev Density Theorem: Every place $w$ of $K$, outside $S(K)$, splits in some $E \in \mathcal{I}_K$, so the union contains $\Gamma_{K_w}$. To check property (b), fix prime degree extensions $K$ and $K'$ as above. Note that
$$
\text{$(\rho_K|_{\Gamma_{KK'}})|_{\Gamma_{EE'}}\simeq (\rho_{K'}|_{\Gamma_{KK'}})|_{\Gamma_{EE'}}$, 
$\y$ $\forall E \in \mathcal{I}_K$, $\y$ $\forall E' \in \mathcal{I}_{K'}$.}
$$
We finish the proof by showing that the union of these $\Gamma_{EE'}$ is dense in $\Gamma_{KK'}$.
Let $w$ be an arbitrary place of $KK'$ such that $w|_F$ does not lie in $S$. Choose an extension $E \in \mathcal{I}_K$ linearly disjoint from $KK'$ over $K$, in which $w|_{K}$ splits. Then pick an 
extension $E' \in \mathcal{I}_{K'}$ linearly disjoint from $EK'$ over $K'$, in which $w|_{K'}$ splits.
By elementary splitting theory, $w$ splits in $EK'$, and any place of $EK'$ above $w$ splits in $EE'$.
Consequently, $w$ splits in $EE'$, see the diagram:
$$
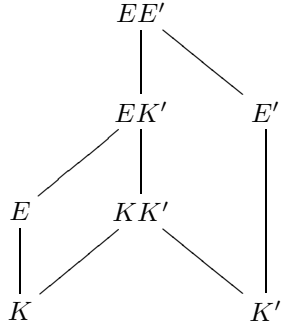
\xymatrix{
& EE' \ar@{-}[d] \ar@{-}[rd]\\
& EK' \ar@{-}[d] \ar@{-}[ld] & E' \ar@{-}[dd]\\
E \ar@{-}[d] & KK' \ar@{-}[ld] \ar@{-}[rd]\\
K & & K'
}
$$
The union then contains the Galois group of $(KK')_w$. Done by Cebotarev. $\square$

\medskip

\noindent The previous result should be compared to Corollary 1.2 in [Har].

\subsection{Galois representations associated to $\pi$}

Let $\pi$ be the globally generic cusp form on $\GSp(4)$ introduced earlier, $\Pi$ its lift to $\GL(4)$, and let $\Pi_E$ be the base change of $\Pi$ to $\GL(4)$ over a CM extension $E$ of $F$. Recall that, for certain algebraic Hecke characters $\chi$ of $E$, the twisted representation $\Pi_E(\chi)$ is conjugate self-dual.
We consider the representations
$$
\rho_{\pi,\iota,E}\overset{\text{df}}{=}\rho_{\Pi_E(\chi),\iota}\otimes \rho_{\check{\chi},\iota}.
$$
Up to isomorphism, this is independent of $\chi$. Indeed, for each place $w\nmid \ell$ of $E$,
$$
\iota\text{WD}(\rho_{\pi,\iota,E}|_{W_{E_w}})^{F-ss}\simeq \text{rec}_4(\Pi_{E,w}\otimes |\det|^{-\frac{3}{2}}).
$$
We only consider CM extensions $E$, in which $v_0$ splits, such that $\Pi_E$ is cuspidal. Here $v_0$ is the place of $F$ where $\pi_{v_0}$ is of Steinberg type. This collection $\mathcal{I}$ is certainly strongly $\varnothing$-general, according to the example in the previous section. Moreover, the family of $4$-dimensional Galois representations $\rho_{\pi,\iota,E}$ satisfies the patching conditions (a) and (b).
For example, to check the {\it{Galois invariance}}, 
$$
\rho_{\pi,\iota,E}^{\theta}\simeq \rho_{\Pi_E(\chi)^{\theta},\iota}\otimes \rho_{\check{\chi}^{\theta},\iota}
\simeq \rho_{\Pi_E(\chi)^{\vee},\iota}\otimes \rho_{\check{\chi}^{\theta},\iota}
\simeq \rho_{\Pi_E(\chi^{\theta}),\iota}\otimes \rho_{\check{\chi}^{\theta},\iota},
$$
by our choice of $\chi$. Taking $\chi^{\theta}$ instead of $\chi$, then shows that $\rho_{\pi,\iota,E}^{\theta}$ is isomorphic to $\rho_{\pi,\iota,E}$ by the aforementioned independence.
Alternatively, one can use the local description of $\rho_{\pi,\iota,E}$ above at the unramified places, and the fact that $\Pi_E$ is a base change from $F$. To check the {\it{compatibility}}, note that for $w \nmid \ell$,
$$
\iota\text{WD}(\rho_{\pi,\iota,E}|_{W_{(EE')_w}})^{F-ss}\simeq \text{rec}_4(\Pi_{v,(EE')_w}\otimes |\det|^{-\frac{3}{2}}),
$$
and similarly for $\rho_{\pi,\iota,E'}$. See Lemma VII.2.6 in [HT]. Now (b) follows from Cebotarev. By the patching lemma, we finally get a continuous representation
$$
\text{$\rho_{\pi,\iota}:\Gamma_F \rightarrow \GL_4(\bar{\Q}_{\ell})$, $\y$ $\rho_{\pi,\iota}|_{\Gamma_E}\otimes \rho_{\chi,\iota}\simeq \rho_{\Pi_E(\chi),\iota}$.}
$$
It is {\it{irreducible}}, since $\rho_{\Pi_E(\chi),\iota}$ is known to be irreducible [TY], and satisfies:
$$
\iota\text{WD}(\rho_{\pi,\iota}|_{W_{F_v}})^{F-ss}\simeq \text{rec}_{\text{GT}}(\pi_v\otimes |c|^{-\frac{3}{2}}),
$$
at each finite place $v \nmid \ell$ of $F$. Here $\text{rec}_{\text{GT}}$ is the local Langlands correspondence for $\GSp(4)$, as defined by Gan and Takeda in [GT]. To see this, pick any $E$ in which $v$ splits, and use the local description of $\rho_{\pi,\iota,E}$ together with the fact that $\Pi$ is a {\it{strong}} lift of $\pi$. From the list of properties of $\rho_{\Pi_E(\chi),\iota}$, we then read off:

\begin{itemize}
\item Let $v\nmid \ell$ be a finite place of $F$, and let $\alpha$ be an eigenvalue of $\rho_{\pi,\iota}(\sigma)$ for some $\sigma \in W_{F_v}$. Then $\alpha$ belongs to $\bar{\Q}$, and for every embedding $\bar{\Q}\hookrightarrow \C$,
$$
|\alpha| \in q_v^{\frac{\Z}{2}}.
$$
\item Let $v\nmid \ell$ be a finite place, with $\pi_v$ {\it{unramified}}, and let $\alpha$ be an eigenvalue of $\rho_{\pi,\iota}(\Frob_v)$. Then $\alpha$ belongs to $\bar{\Q}$, and for every embedding $\bar{\Q}\hookrightarrow \C$,
$$
|\alpha|=q_v^{\frac{w+3}{2}}.
$$
\item The representation $\rho_{\pi,\iota}$ is potentially semistable at any finite place $v|\ell$. Moreover, $\rho_{\pi,\iota}$ is crystalline at a finite place $v|\ell$ when $\pi_v$ is unramified.
\end{itemize}
For the second part, we recall that $\chi$ is an algebraic Hecke character with infinity types $z^a\bar{z}^b$, where $a+b=w$. In particular, for the {\it{unitary}} twist $\pi^{\circ}$ we have:
$$
L_v(s-\frac{1}{2}(w+3),\pi^{\circ},\text{spin})=\det(1-\iota\rho_{\pi,\iota}(\text{Frob}_v)\cdot q_v^{-s})^{-1}
$$
at all places $v \nmid \ell$ where $\pi^{\circ}$ is unramified. Note that, by twisting $\rho_{\pi,\iota}$ with integral powers of the cyclotomic character $\chi_{\text{cyc}}$, we may alter the motivic weight $w+3$ by any even integer. We compare with the motivic weight $k_1+k_2-3$ in [Wei].

\medskip

\noindent {\it{Temperedness of $\pi^{\circ}$:}} From the above, it follows immediately that $\pi^{\circ}$ has unitary Satake parameters at all places $v \nmid \ell$ where $\pi^{\circ}$ is unramified. In fact, $\pi^{\circ}$ is tempered at {\it{every}} place $v$: Indeed, by Corollary VII.1.11 in [HT], we know that $\Pi_E$ is essentially tempered everywhere. That is, $\phi_{\pi_v}|_{W_{E_w}}$ has bounded image in $\GL_4(\C)$ for every finite place $w$ of $E$. Consequently, the same holds for $\phi_{\pi_v}|_{W_{F_v}}$.

\medskip

\noindent {\it{The image of $\rho_{\pi,\iota}$:}} Since the eigenvalues of $\rho_{\pi,\iota}(\text{Frob}_v)$ coincide with the integral Satake parameters of $\pi_v$, for finite $v\nmid \ell$ where $\pi_v$ is unramified, Cebotarev yields:
$$
\text{$\rho_{\pi,\iota}^{\vee}\simeq \rho_{\pi,\iota}\otimes \chi^{-1}$, $\y$ $\chi\overset{\text{df}}{=}\omega_{\pi^{\circ}}\cdot\chi_{\text{cyc}}^{-w-3}$,}
$$
where we confuse $\omega_{\pi^{\circ}}$ with its finite order $\ell$-adic avatar. In other words, the space of $\rho_{\pi,\iota}$ has a non-degenerate bilinear form preserved by $\Gamma_F$ with similitude $\chi$. We have already observed that $\rho_{\pi,\iota}$ is irreducible, so by Schur's lemma this bilinear form must be symmetric {\it{or}} symplectic. Thus, the image of $\rho_{\pi,\iota}$ can always be conjugated into $\GO(4)$ or $\GSp(4)$. Under our running assumptions on $\pi$, in fact into the latter: Otherwise, by local-global compatibility at the place $v_0\nmid \ell$, the $L$-parameter of $\Pi_{v_0}$ is of orthogonal type. That is,
it maps
$$
\text{rec}(\Pi_{v_0}): W_{F_{v_0}}'=W_{F_{v_0}}\times \SL(2,\C) \rightarrow \GO(4,\C).
$$ 
However, $\Pi_{v_0}$ is the transfer of $\pi_{v_0}$, so $\text{rec}(\Pi_{v_0})$ also preserves a symplectic form on $\C^4$. Now, $\Pi_{v_0}$ is a generalized Steinberg representation, and one verifies that
$$
Z_{\GL(4,\C)}(\text{im}(\text{rec}(\Pi_{v_0})))=\C^*
$$
by an easy computation. Indeed, $\text{rec}(\Pi_{v_0})$ is of the form $\phi \boxtimes S_d$, where $\phi$ is an irreducible representation of $W_{F_{v_0}}$, and $S_d$ is the $d$-dimensional irreducible representation of $\SL(2,\C)$. Ergo, the above symplectic form must agree with the orthogonal form up to a scalar. This is a contradiction. The symplecticity of $\rho_{\pi,\iota}$, just shown, is a special case of Theorem F on p. 6 in [CCl] when $\omega_{\pi}$ is trivial. This result from [CCl] has recently been generalized to the CM case in [BCh]. When $F=\Q$, the symplecticity of $\rho_{\pi,\iota}$ is shown in [Wei], for globally generic $\pi$, using Poincare duality. Indeed, by [Sou], $\pi$ occurs with multiplicity one, so $\rho_{\pi,\iota}$ can be realized as the $\pi_f$-isotypic component of $H^3$ of a Siegel threefold. The cup product pairing then provides the desired symplectic form. 

\medskip

\noindent {\it{Baire category theory:}} To check that the image of $\rho_{\pi,\iota}$ is in fact contained in $\GSp_4(L)$, for some finite extension $L$ over $\Q_{\ell}$, we invoke the Baire Category Theorem: Every locally compact Hausdorff space is a Baire space (that is, the union of any countable collection of closed sets with empty interior has empty interior). We will apply it to the compact subgroup $\rho_{\pi,\iota}(\Gamma_F)$ inside $\GSp_4(\bar{\Q}_{\ell})$.
$$
\rho_{\pi,\iota}(\Gamma_F)={\bigcup}_{\text{$L/\Q_{\ell}$ finite}} \rho_{\pi,\iota}(\Gamma_F) \cap \GSp_4(L)
$$
is a countable union of closed subgroups, since each $L$ is complete. Therefore, 
$$
\text{$\rho_{\pi,\iota}(\Gamma_F) \cap \GSp_4(L)$ has {\it{non-empty}} interior,}
$$
for some $L$, and hence this is an open subgroup. That is, the image of $\Gamma_M$ for some finite extension $M$ over $F$. In particular, it has {\it{finite}} index in $\rho_{\pi,\iota}(\Gamma_F)$. By enlarging $L$, to accomodate the finitely many coset representatives, we can arrange for the image of $\rho_{\pi,\iota}$ to be contained in the $L$-rational points $\GSp_4(L)$.

\medskip

\noindent {\it{Total oddness:}} $\chi(c)=-1$ for every complex conjugation $c \in \Gamma_F$ from $\bar{\Q}\hookrightarrow \C$.

\medskip

\noindent {\it{Hodge-Tate weights:}} Let us fix an embedding $s: F \rightarrow \bar{\Q}_{\ell}$, and let $v$ be the associated place of $F$ above $\ell$. We wish to compute the Hodge-Tate weights of $\rho_{\pi,\iota}$ restricted to $\Gamma_{F_v}$, where $F_v=s(F)^-$. That is, for each integer $j$, evaluate
$$
\dim_{\bar{\Q}_{\ell}}\text{gr}^j(\rho_{\pi,\iota}\otimes_{F_v}B_{dR})^{\Gamma_{F_v}}.
$$
We will reduce this to the analogous result for $\rho_{\Pi_E(\chi),\iota}$ already mentioned. Thus, we fix a CM extension $E$, in which $v$ splits. Once and for all, we fix a divisor $w$ of $v$, and look at a corresponding embedding $\tilde{s}:E \rightarrow \bar{\Q}_{\ell}$ over $s$. This canonically identifies $E_w=\tilde{s}(E)^-$ with $F_v$. Now note that, for characters $\chi$ as above,
$$
\rho_{\pi,\iota}|_{\Gamma_{F_v}}\otimes \rho_{\chi,\iota}|_{\Gamma_{E_w}}\simeq \rho_{\Pi_E(\chi),\iota}|_{\Gamma_{E_w}}.
$$
Therefore, we first record the Hodge-Tate weight of $\rho_{\chi,\iota}|_{\Gamma_{E_w}}$. The associated complex embedding $\iota(\tilde{s})$ defines an infinite place of $E$, where $\chi$ has the form $z^a\bar{z}^b$. As is well documented elsewhere in the literature, for example in [Bla], the Hodge-Tate weight is then $-b$, with our choice of normalization. Therefore,
$$
\dim_{\bar{\Q}_{\ell}}\text{gr}^j(\rho_{\pi,\iota}\otimes_{F_v}B_{dR})^{\Gamma_{F_v}}=
\dim_{\bar{\Q}_{\ell}}\text{gr}^{j-b}(\rho_{\Pi_E(\chi),\iota}\otimes_{E_w}B_{dR})^{\Gamma_{E_w}},
$$
since $D_{dR}$ is a $\otimes$-functor. It remains to find the highest weight of the $\mathcal{V}$ with
$$
H^{\bullet}(\frak{g},K; \Pi_E(\chi)_{\infty}\otimes \mathcal{V}^*)\neq 0.
$$
More precisely, we let $\sigma=\iota(s) \in \Sigma$, and consider the two complex embeddings 
$\{\tilde{\sigma},\tilde{\sigma}^c\}$ of $E$ extending $\sigma$. Here $\tilde{\sigma}=\iota(\tilde{s})$.
In our earlier notation, we need to compute the quadruple $\mu(\tilde{\sigma}^c)$. For this, we follow the proof of Lemma 3.14 on p. 114 in [Cl2]: We consider the local component of $\Pi_E(\chi)$ at the infinite place of $E$ above $\sigma$. We know its $L$-parameter, so according to p. 113 in [Cl2]:
$$
\begin{cases}
\mu_1(\tilde{\sigma}^c)=b+3-\frac{1}{2}({\bf{w}}+n'), \\
\mu_2(\tilde{\sigma}^c)=b+2-\frac{1}{2}({\bf{w}}+n), \\
\mu_3(\tilde{\sigma}^c)=b+1-\frac{1}{2}({\bf{w}}-n), \\
\mu_4(\tilde{\sigma}^c)=b+0-\frac{1}{2}({\bf{w}}-n'). \\
\end{cases}
$$
Here we have introduced $n=\nu_1-\nu_2$ and $n'=\nu_1+\nu_2$. Moreover, the motivic weight $w+3$ is denoted by ${\bf{w}}$. From the above, and the result from section 4.2, we deduce that the Hodge-Tate weights of $\rho_{\pi,\iota}|_{\Gamma_{F_v}}$ are given by the sequence:
$$
\frac{1}{2}({\bf{w}}-n')<\frac{1}{2}({\bf{w}}-n)<\frac{1}{2}({\bf{w}}+n)<\frac{1}{2}({\bf{w}}+n').
$$
In particular, they are distinct. We will rewrite this slightly. For each $\sigma \in \Sigma$, 
$$
\delta=\delta(\sigma)\overset{\text{df}}{=}\frac{1}{2}({\bf{w}}-n')=\frac{1}{2}(w-\mu_1-\mu_2)\in \Z.
$$
With this notation, the set of Hodge-Tate weights takes the following form:
$$
\text{HT}(\rho_{\pi,\iota}|_{\Gamma_{F_v}})=\{\delta, \nu_2+\delta, \nu_1+\delta, \nu_1+\nu_2+\delta\}.
$$
In the case $F=\Q$ it is customary to take ${\bf{w}}=k_1+k_2-3$, that is, $\delta=0$. In this case, we recover the Hodge types given in Theorem III on p. 2 in [Wei].

\subsection{Consequences of local-global compatibility}

\noindent {\it{Parahoric subgroups:}} We fix a finite place $v$ of $F$, and define certain compact open subgroups of $\GSp_4(F_v)$, known as the parahoric subgroups. They arise as stabilizers of points in the Bruhat-Tits builidng. We refer to [Tit] for a general discussion. First, we have the {\it{hyperspecial}} maximal compact subgroup
$$
K=K_v\overset{\text{df}}{=}\GSp_4(\mathcal{O}_v).
$$
Inside of it, we have the pullbacks of the two parabolics via the reduction map:
$$
\begin{cases}
J_P=J_{P,v}\overset{\text{df}}{=}\{\text{$k \in K$: $k$ (mod $v$) $\in P(\F_v)$}\}, \\
J_Q=J_{Q,v}\overset{\text{df}}{=}\{\text{$k \in K$: $k$ (mod $v$) $\in Q(\F_v)$}\}. \\
\end{cases}
$$
They are usually called the {\it{Siegel}} and {\it{Klingen}} parahoric, respectively. Moreover,
$$
I=I_v\overset{\text{df}}{=}J_P\cap J_Q=\{\text{$k \in K$: $k$ (mod $v$) $\in B(\F_v)$}\}
$$
is called the (upper-triangular) {\it{Iwahori}} subgroup. Furthermore, let us introduce 
$$
\text{$\eta=\eta_v\overset{\text{df}}{=}\begin{pmatrix} & & 1 & \\ & & & 1 \\ \varpi_v & & & \\ & \varpi_v & & \end{pmatrix}$, $\y$ $c(\eta)=-\varpi_v$, $\y$ $\eta^2=\varpi_v\cdot I_4$.}
$$
This $\eta$ is occasionally referred to as the {\it{Atkin-Lehner}} element. It depends on the choice of a uniformizer $\varpi_v$, but it is well-defined modulo the maximal compact subgroup of the torus $T_c$. 
An explicit calculation shows that the $\eta$-conjugate $\eta K \eta^{-1}$ is another hyperspecial maximal compact subgroup containing the Siegel parahoric. There is one more maximal compact subgroup containing $I$. Namely,
$$
\tilde{K}=\tilde{K}_v\overset{\text{df}}{=}\text{subgroup generated by $J_Q$ and its $\eta$-conjugate.}
$$
This is the {\it{paramodular}}, or non-special, subgroup. We have thus described all parahoric subgroups of $\GSp_4(F_v)$ up to conjugacy. Note that $\eta$ normalizes $\tilde{K}$ and $J_P$, hence $I$. Therefore, if $J$ is any one of these subgroups, and $\pi$ is an irreducible admissible representation, $\pi(\eta)$ is a well-defined operator on the $J$-invariants $\pi^J$. Its square is multplication by $\omega_{\pi}(\varpi_v)$. In particular, if $\omega_{\pi}=1$, the eigenvalues of $\pi(\eta)$ are $\pm 1$. 
These are tabulated in table 3 on p. 16 in [Sch].

\medskip

\noindent {\it{Nilpotent orbits:}} Suppose $\pi$ is an irreducible representation of $\GSp_4(F_v)$, assumed to be Iwahori-spherical. That is, $\pi^I \neq 0$. Then its $L$-parameter is given by a semisimple element in $\GSp_4(\C)$ plus a compatible nilpotent element $N$ in the Lie algebra $\frak{gsp}_4(\C)$, both viewed up to conjugacy. Here we list the finitely many possibilities for $N$. First, by the theory of the Jordan normal form,
$$
\text{$\begin{pmatrix}0 & 1 & & \\ & 0 & & \\ & & 0 & \\ & & & 0\end{pmatrix}$,
$\begin{pmatrix}0 & 1 & & \\ & 0 & & \\ & & 0 & 1\\ & & & 0\end{pmatrix}$,
$\begin{pmatrix}0 & 1 & & \\ & 0 & 1 & \\ & & 0 & \\ & & & 0\end{pmatrix}$,
$\begin{pmatrix}0 & 1 & & \\ & 0 & 1 & \\ & & 0 & 1\\ & & & 0\end{pmatrix}$}
$$
represent the non-trivial nilpotent classes in $\frak{gl}_4(\C)$. An explicit computation verifies that the third representative {\it{cannot}} be conjugated into $\frak{gsp}_4(\C)$. Recall that $\frak{sp}_4(\C)$ consists of $X$ such that $JX$ is symmetric. Here $J$ is the anti-diagonal symplectic form fixed throughout the paper.
However, the other three representatives {\it{are}} in fact symplectic in this sense: They are conjugate to
$$
\text{$\mathcal{N}_1\overset{\text{df}}{=}\begin{pmatrix}0 & & & \\ & 0 & 1 & \\ & & 0 & \\ & & & 0\end{pmatrix}$,
$\mathcal{N}_2\overset{\text{df}}{=}\begin{pmatrix}0 & 1 & & \\ & 0 &  & \\ & & 0 & -1\\ & & & 0\end{pmatrix}$,
$\mathcal{N}_3\overset{\text{df}}{=}\begin{pmatrix}0 & 1 & & \\ & 0 & 1 & \\ & & 0 & -1\\ & & & 0\end{pmatrix}$}
$$
respectively. Again, this follows immediately from the theory of normal forms. Each $\mathcal{N}_i$ is contained in $\frak{sp}_4(\C)$. Note that $\mathcal{N}_i$ has rank $i$. It may be useful to observe that 
$\mathcal{N}_1$ is a root vector for the {\it{long}} simple root $\beta(t)=2t_2$, whereas $\mathcal{N}_2$ is a root vector for the {\it{short}} simple root $\alpha(t)=t_1-t_2$. Their sum is $\mathcal{N}_3$. By the Jacobson-Morozov theorem, [Jac, Theorem 3], the $\mathcal{N}_i$ in fact represent the $\GSp_4(\C)$-conjugacy classes of nilpotent elements in $\frak{gsp}_4(\C)$. To aid comparison with [Sch], let us make the following remark: For the $L$-parameter of $\pi$ to have a more transparent semisimple part, one often takes a different set of representatives, see p. 6 in [Sch]. For example, the two nilpotent elements
$$
\text{$\mathcal{N}_1'\overset{\text{df}}{=}\begin{pmatrix}0 & & & 1\\ & 0 &  & \\ & & 0 & \\ & & & 0\end{pmatrix}$, $\y$ 
$\mathcal{N}_2'\overset{\text{df}}{=}\begin{pmatrix}0 & & & 1\\ & 0 & 1 & \\ & & 0 & \\ & & & 0\end{pmatrix}$,}
$$
are $\GSp_4(\C)$-conjugate to $\mathcal{N}_1$ and $\mathcal{N}_2$ respectively. This is easy to check.

\medskip

\noindent {\it{Iwahori-spherical generic representations:}} By a well-known result of Casselman, see [Car] for a nice review, the Iwahori-spherical $\pi$ are precisely the constituents of unramified principal series.
The way they decompose, in the case of $\GSp_4$, was written out explicitly in [ST]. One gets $17$ families of representations,
$$
I, \y IIa, \y IIb, \y IIIa, \y IIIb, \y IVa-IVd, \y Va-Vd, \y VIa-VId,
$$
according to how the unramified principal series breaks up. We refer to pages 6-8 in [Sch] for a precise definition of each class of representations. Fortunately, here we are only interested in the {\it{generic}} representations. That is, the $6$ classes
$$
I, \y IIa, \y IIIa, \y IVa, \y Va, \y VIa.
$$
We will briefly recall the definition of each of these classes. To do that, we will first introduce the notation used in [ST]: For three quasi-characters $\chi_i$, let
$$
\chi_1 \times \chi_2 \rtimes \chi_3\overset{\text{df}}{=}\text{Ind}_B(\chi_1\otimes\chi_2\otimes \chi_3)
$$
be the principal series for $\GSp_4$ obtained by normalized induction from 
$$
\chi_1\otimes\chi_2\otimes \chi_3: t \mapsto \chi_1(t_1)\chi_2(t_2)\chi_3(c(t)).
$$
Similarly, if $\tau$ is an irreducible representation of $\GL_2$, and $\chi$ is a character, let
$$
\text{$\tau \rtimes \chi\overset{\text{df}}{=}\text{Ind}_P(\tau \otimes \chi)$, $\y$ 
$\chi \rtimes \tau \overset{\text{df}}{=}\text{Ind}_Q(\chi \otimes \tau)$.}
$$
Again, we use {\it{normalized}} induction, and we identify the Levi subgroups of $P$ and $Q$ with $\GL_2\times \GL_1$ in the natural way. Having this notation at hand, we describe the generic classes of Iwahori-spherical representations discussed above:

\begin{enumerate}
\item[($I$)] $\pi=\chi_1 \times \chi_2 \rtimes \chi_3$,
\item[($IIa$)] $\pi=\text{St}_{\GL(2)}(\chi_1)\rtimes \chi_2$,
\item[($IIIa$)] $\pi=\chi_1 \rtimes \text{St}_{\GL(2)}(\chi_2)$,
\item[($IVa$)] $\pi=\text{St}_{\GSp(4)}(\chi)$, 
\item[($Va$)] $\pi\subset \text{St}_{\GL(2)}(\nu^{1/2}\xi_0)\rtimes \nu^{-1/2}\chi$
\item[($VIa$)] $\pi \subset {\bf{1}}\rtimes \text{St}_{\GL(2)}(\chi)$
\end{enumerate}

\noindent Here all the characters $\chi$ and $\chi_i$ are unramified. Moreover, $\nu$ denotes the normalized absolute value, and $\xi_0$ is the non-trivial unramified quadratic character. According to Table 1 on p. 9 in [Sch], only $IVa$ and $Va$ are discrete series. Type $VIa$ representations are the analogues of {\it{limits}} of discrete series. Table 3 on p. 16 in [Sch] lists the dimensions of the parahoric fixed spaces for all 17 families above, plus additional data such as the Atkin-Lehner eigenvalues when $\omega_{\pi}=1$. Below, we concatenate parts of Table 2 and parts of Table 3 in [Sch]. That is,

{\center{

\begin{tabular}{|c|c|c|c|c|c|c|c|c|}
\hline
  type  & $N$ & $K$ & $\tilde{K}$ & $J_P$ & $J_Q$ & $I$ \\
 \hline \hline
  $I$     &  0 & 1 &  2 & 4 & 4 & 8\\
\hline
  $IIa$ 
   &  $\mathcal{N}_1$ & 0 & 1 & 1 & 2 & 4\\
 \hline
  $IIIa$  
  &  $\mathcal{N}_2$ &  0 & 0 & 2 & 1& 4\\
 \hline
$IVa$   & $\mathcal{N}_3$ & 0 &  0 & 0 &0 & 1\\
 \hline
$Va$   & $\mathcal{N}_2$ & 0 & 0 & 0 & 1 & 2\\
 \hline
$VIa$   & $\mathcal{N}_2$ & 0 & 0 &1  & 1& 3\\
 \hline
\end{tabular}

{\center{Table A: Parahoric fixed spaces and monodromy}}

}}

\medskip

\noindent We note that the assignment of monodromy operators in [Sch] is compatible with the local Langlands correspondence, as defined by Gan and Takeda via theta correspondence. This follows from the explicit calculations in section 12 of [GT], see their remarks on p. 33. By local-global compatibility, we deduce: 

\begin{cor}
Let $\rho_{\pi,\iota}$ be the Galois representation attached to a globally generic cusp form $\pi$ as above. Let $v \nmid \ell$ be a finite place of $F$ such that $\pi_v$ is Iwahori-spherical and ramified. Then $\rho_{\pi,\iota}|_{I_{F_v}}$ acts unipotently. In fact, the image $\rho_{\pi,\iota}(I_{F_v})$ is topologically generated by a $\GSp_4(\bar{\Q}_{\ell})$-conjugate of $\exp(\mathcal{N}_i)$, where
$$
i=\begin{cases}
1, \y \text{$\pi_v$ of type IIa,} \\
2, \y \text{$\pi_v$ of type IIIa, Va, or VIa,} \\
3, \y \text{$\pi_v$ of type IVa.} \\
\end{cases}
$$
In particular, we have the following consequences conjectured in [GTi] and [SU]:

\begin{itemize}
\item $\pi_v$ of Steinberg type $\Longleftrightarrow$ monodromy has rank $3$. 
\item $\pi_v$ has a unique $J_Q$-fixed line $\Longleftrightarrow$ monodromy has rank $2$.
\item $\pi_v$ para-spherical $\Longleftrightarrow$ monodromy has rank $1$.
\end{itemize}

\end{cor}

\noindent $Proof$. This follows immediately from Table A above. $\square$

\medskip

\noindent The first two consequences are part of the Conjecture on p. 11 in [GTi]. Note that part 4 of that Conjecture is false: If $\pi_v$ has a unique $J_P$-fixed line, one cannot deduce that monodromy has rank one. The last consequence is Conjecture 3.1.7 on p. 41 in [SU], for globally generic $\pi$ as above. Skinner and Urban used the holomorphic analogue as a substitute for deep results of Kato, in order to study Selmer groups for certain modular forms of square-free level. Another application to the Bloch-Kato conjecture in this context, contingent on the holomorphic analogue of the third consequence above, was given in [Sor].

\medskip

\noindent {\it{Supercuspidal generic representations:}} According to Table 2 on p. 51 in [GT], there are two types of supercuspidal generic representations $\pi$ of $\GSp_4$. Firstly,
$$
\pi=\theta((\sigma\otimes \sigma')^+),
$$
for distinct supercuspidal representations $\sigma \neq \sigma'$ on $\GL_2$. In this case, the lift to $\GL_4$ is the isobaric sum $\sigma \boxplus \sigma'$. Secondly, if $\pi$ is not a lift from $\GO_{2,2}$, when lifted to $\GL_4$ it remains supercuspidal. Again, by local-global compatibility:

\begin{cor}
Let $\rho_{\pi,\iota}$ be the Galois representation attached to a globally generic cusp form $\pi$ as above. Let $v \nmid \ell$ be a finite place of $F$ such that $\pi_v$ is supercuspidal. Then $\rho_{\pi,\iota}$ is trivial on some finite index subgroup of $I_{F_v}$. Moreover,
$$
\text{$\pi_v$ is \underline{not} a lift from $\GO_{2,2}$ $\Longleftrightarrow$ $\rho_{\pi,\iota}|_{W_{F_v}}$ is irreducible.}
$$
On the contrary, when $\pi_v$ \underline{is} a lift from $\GO_{2,2}$, the restriction $\rho_{\pi,\iota}|_{W_{F_v}}^{\text{ss}}$ breaks up as a sum of two non-isomorphic irreducible two-dimensional representations.
\end{cor}

\noindent $Proof$. This follows from the foregoing discussion. $\square$

\medskip

\noindent For a moment, let us continue with the setup of the previous Corollary. The exponent of the Artin conductor of $\rho_{\pi,\iota}|_{I_{F_v}}$ is defined by the standard formula:
$$
\frak{f}(\rho_{\pi,\iota}|_{I_{F_v}})\overset{\text{df}}{=}\sum_{i=0}^{\infty}\frac{1}{[\tilde{I}_{F_v}:\tilde{I}_{F_v,i}]}\cdot \text{codim}_{\bar{\Q}_{\ell}}(\rho_{\pi,\iota}^{\tilde{I}_{F_v,i}})\in \Z,
$$
where $\tilde{I}_{F_v,i}$ is the $i$th ramification group in $\tilde{I}_{F_v}$, in turn some finite quotient of $I_{F_v}$ through which $\rho_{\pi,\iota}$ factors. This sum is finite. The exponent of the Swan conductor, $\frak{f}_{\text{Swan}}$, is given by the same formula except the summation starts at $i=1$. By irreducibility, it is easy to see that $\rho_{\pi,\iota}$ has no nonzero $I_{F_v}$-invariants:
$$
\frak{f}_{\text{Swan}}(\rho_{\pi,\iota}|_{I_{F_v}})=\frak{f}(\rho_{\pi,\iota}|_{I_{F_v}})-4.
$$
We wish to relate this to the {\it{depth}} of $\pi_v$, a non-negative rational number measuring its wild ramification. We very briefly recall the definition: Let $G$ be (the rational points of) a connected reductive group over $F_v$, and let $x$ be a point on its extended Bruhat-Tits building [Tit]. Its stabilizer $G_x$ is the corresponding parahoric subgroup. In [MP], Moy and G. Prasad defined an exhaustive descending filtration of $G_x$, consisting of open subgroups $G_{x,r}$ parametrized by non-negative real numbers $r$. They then defined a pro-$p$ group $G_{x,r^+}$ to be the union of the $G_{x,s}$ for $s>r$. The breaks $r$, where $G_{x,r^+}$ is a proper subgroup of $G_{x,r}$, is known to form an unbounded discrete subset of $\R$. The depth of $\pi$ is
$$
\text{depth}(\pi)\overset{\text{df}}{=}\text{inf}\{\text{$r$: $\pi^{G_{x,r^+}}\neq 0$, some $x$}\}\in \Q.
$$
Here $\pi$ is any irreducible admissible representation of $G$. Our goal is to show:

\begin{prop}
Let $\rho_{\pi,\iota}$ be the Galois representation attached to a globally generic cusp form $\pi$ as above. Let $v \nmid \ell$ be a finite place of $F$ such that $\pi_v$ is supercuspidal, and \underline{not} a lift from $\GO_{2,2}$. Then we have the following identity:
$$
\frak{f}_{\text{Swan}}(\rho_{\pi,\iota}|_{I_{F_v}})=4 \cdot \text{depth}(\pi_v).
$$
\end{prop}

\noindent $Proof$. As is well-known, see [Tat], the Artin conductor fits into the $\epsilon$-factor:  
$$
\epsilon(s,\iota\text{WD}(\rho_{\pi,\iota}|_{W_{F_v}}),\psi)=\epsilon(0,\iota\text{WD}(\rho_{\pi,\iota}|_{W_{F_v}}),\psi)\cdot q_v^{-s(\frak{f}(\rho_{\pi,\iota}|_{I_{F_v}})+4n(\psi))}.
$$
Here $\psi$ is some fixed non-trivial character of $F_v$, and $n(\psi)$ is the largest $n$ such that $\psi$ is trivial on $\frak{p}_v^{-n}$. Similarly, if $\Pi_v$ is the supercuspidal lift of $\pi_v$ to $\GL_4$,  
$$
\epsilon(s,\Pi_v,\psi)=\epsilon(0,\Pi_v,\psi)\cdot q_v^{-s(\frak{f}(\Pi_v)+4n(\psi))}.
$$
Here $\frak{f}(\Pi_v)$ is the standard conductor of $\Pi_v$, that is, the smallest $f$ such that $\Pi_v$ has nonzero vectors fixed by the subgroup consisting of elements in $\GL_4(\mathcal{O}_v)$ whose last row is congruent to $(0,\ldots, 0,1)$ modulo $\frak{p}_v^f$. Hence, we deduce that
$$
\frak{f}(\rho_{\pi,\iota}|_{I_{F_v}})=\frak{f}(\Pi_v),
$$
since the local Langlands correspondence for $\GL_4$ preserves $\epsilon$-factors. The determinant twist can be ignored. Now, the key ingredient is the following formula due to Bushnell and Frolich [BF], which holds for supercuspidals on any $\GL_n$,
$$
\text{$\frak{f}(\Pi_v)=n \cdot \text{depth}(\Pi_v)+n$, $\y$ $n=4$.}
$$
We note, in passing, that this formula was generalized to the square integrable case in [LR]. It remains to explain why $\pi_v$ and $\Pi_v$ have the same depth. Keep in mind that $\Pi_v\otimes \omega_{\pi_v}$ is the theta lift of $\pi_v$ to $\GSO_{3,3}$. Now invoke the main result from [Pan], suitably extended to incorporate similitudes. For this last step, Lemma 2.2 on p. 7 in [GT] is very useful. We omit the details. $\square$

\medskip

\noindent {\it{Tame ramification:}} In the previous Proposition, let us take $\pi_v$ to have depth zero. In other words, $\pi_v$ has nonzero vectors fixed by the pro-unipotent radical of some parahoric subgroup. In this special case, $\rho_{\pi,\iota}|_{I_{F_v}}$ factors through the tame quotient $I_{F_v}^t$, a pro-cyclic group of pro-order prime-to-$p$. More concretely, $\rho_{\pi,\iota}|_{I_{F_v}}$ is the direct sum of four $\ell$-adic characters of $I_{F_v}$, of finite order prime-to-$p$. Below, we will state a related result, due to Genestier and Tilouine in the rational case. Let $\chi$ be any complex character of $\F_v^*$ and view it as a character 
$$
\chi: J_{Q,v} \rightarrow \GL(2,\F_v)\times \F_v^* \rightarrow \F_v^* \rightarrow \C^*,
$$
where the second map is projection. If $\pi$ is an irreducible admissible representation of $\GSp_4$, we will be looking at the space $\pi^{J_Q,\chi}$ of vectors on which $J_Q$ acts via the character $\chi$.
When $\chi \neq 1$, this is non-trivial only for principal series:

\begin{lem}
Let $\pi$ be an irreducible generic representation of $\GSp_4(F_v)$, and let $\chi$ be a non-trivial character of $\F_v^*$, viewed as a character of $J_Q$ as above. Then $\pi^{J_Q,\chi}$ is nonzero if and only if $\pi$ is a tamely ramified principal series of the form
$$
\pi=\tilde{\chi} \times (\text{unram.}) \rtimes (\text{unram.}),
$$
for some extension $\tilde{\chi}$ of $\chi$ inflated to a tamely ramified character of $\mathcal{O}_v^*$.
\end{lem}

\noindent $Proof$. First, let us assume $\pi^{J_Q,\chi}$ contains nonzero vectors. On such vectors, the Iwahori subgroup $I$ acts via the character $\chi \otimes 1 \otimes 1$. By Roche's construction of types for principal series [Roc], see the formulation on p. 10 in [So2], we deduce that $\pi$ must be a subquotient of a principal series representation of the form
$$
\tilde{\chi} \times \chi_1 \rtimes \chi_2
$$
as in the Lemma. That is, both $\chi_i$ are unramified, and $\tilde{\chi}$ extends $\chi$. Our goal is to show that this principal series is necessarily irreducible. For that, we use the criterion from Theorem 7.9 in [Tad]. Since $\tilde{\chi}$ is ramified, the only way it could be reducible is if $\chi_1=\nu^{\pm1}$. Recall that $\nu$ denotes the normalized absolute value on $F_v$. For simplicity, let us assume that $\chi_1=\nu$. The other case is taken care of by taking the dual. Then, by Lemmas 3.4 and 3.9 in [ST], there is a sequence   
$$
0 \rightarrow \tilde{\chi}\rtimes \text{St}_{\GL(2)}(\sigma) \rightarrow \tilde{\chi} \times \nu \rtimes \nu^{-1/2}\sigma \rightarrow \tilde{\chi}\rtimes {\bf{1}}_{\GL(2)}(\sigma)\rightarrow 0,
$$
where we write $\chi_2$ as $\nu^{-1/2}\sigma$. Both constituents are irreducible. However, the quotient is non-generic, so $\pi$ must be the subrepresentation. It remains to check
$$
(\tilde{\chi}\rtimes \text{St}_{\GL(2)}(\sigma))^{J_Q,\chi}=0.
$$
This is done by explicit calculation: We fix a complete set of representatives,
$$
\text{$Q\backslash G/J_Q=\{1,s_1,s_1s_2s_1\}$, $\y$ $s_1\overset{\text{df}}{=}\begin{pmatrix} & 1 & & \\ 1 & & & \\ & & & 1 \\ & & 1 & \end{pmatrix}$, $\y$ $s_2\overset{\text{df}}{=}\begin{pmatrix} 1 & & & \\ & & 1 & \\ & -1 & & \\ & & & 1\end{pmatrix}$.}
$$
We are inducing from the Klingen parabolic, and an easy argument shows that the $\chi$-vectors in the induced representation are given by the $\chi$-vectors in the representation of the Levi subgroup. In our case, we get the three contributions:  

\begin{itemize}
\item $(\tilde{\chi}\otimes \text{St}_{\GL(2)}(\sigma))^{Q \cap J_Q,\chi}=\tilde{\chi}^{\mathcal{O}^*,\chi}\otimes \text{St}_{\GL(2)}(\sigma)^{\GL_2(\mathcal{O})}=0$,
\item $(\tilde{\chi}\otimes \text{St}_{\GL(2)}(\sigma))^{Q \cap s_1J_Qs_1,\chi}=\tilde{\chi}^{\mathcal{O}^*}\otimes \text{St}_{\GL(2)}(\sigma)^{I,\chi}=0$,
\item $(\tilde{\chi}\otimes \text{St}_{\GL(2)}(\sigma))^{Q \cap s_1s_2s_1J_Qs_1s_2^{-1}s_1,\chi}=\tilde{\chi}^{\mathcal{O}^*,\chi'}\otimes \text{St}_{\GL(2)}(\sigma)^{\GL_2(\mathcal{O})}=0.$
\end{itemize}
Here $\chi'$ is some irrelevant character. This proves the only if part of the Lemma. The converse is easier.
Indeed, by the same observation, it suffices to check that
$$
(\tilde{\chi}\otimes \chi_1\otimes \chi_2)^{B\cap J_Q,\chi}=\tilde{\chi}^{\mathcal{O}^*,\chi}\otimes \chi_1^{\mathcal{O}^*}\otimes \chi_2^{\mathcal{O}^*}\neq 0,
$$
as follows from our assumptions on these characters. This finishes the proof. $\square$

\medskip

\noindent As a last application of local-global compatibility, in conjunction with the previous Lemma, we obtain the following result due to Genestier and Tilouine in the rational case $F=\Q$; compare to the second part of Theorem 2.2.5 in [GTi].  

\begin{cor}
Let $\rho_{\pi,\iota}$ be the Galois representation attached to a globally generic cusp form $\pi$ as above. Let $v \nmid \ell$ be a finite place of $F$ such that $\pi_v^{J_Q,\chi}$ is \underline{nonzero} for some \underline{non-trivial} tamely ramified character $\chi$ of $\mathcal{O}_v^*$. It follows that
$$
\rho_{\pi,\iota}|_{I_{F_v}}=1 \oplus 1 \oplus \chi \oplus \chi.
$$
Here $\chi$ is the character of $I_{F_v}$ obtained via local class field theory. Moreover, one can arrange for the two eigenspaces, for $1$ and $\chi$, to be totally isotropic.
\end{cor}

\noindent $Proof$. The previous Lemma. See also (vi) of Proposition 12.15 in [GT]. $\square$

\medskip

\noindent In [GTi], this result was proved by a careful study of the bad reduction of Siegel threefolds with Klingen level structure at $v$. The eigenspace polarization comes from the cohomology of each irreducible component of the special fiber.

\noindent {\sc{Department of Mathematics, Princeton University, USA.}}

\noindent {\it{E-mail address}}: {\texttt{claus@princeton.edu}}

\end{document}